\documentclass[11pt,a4paper]{article}
\usepackage[T1]{fontenc}
\usepackage{lmodern}
\usepackage[margin=27mm,top=25mm,bottom=26mm,headheight=15pt,headsep=7mm]{geometry}
\usepackage{amsmath,amssymb,amsthm,mathtools}
\usepackage{booktabs,tabularx,array,enumitem}
\usepackage{tikz}
\usetikzlibrary{arrows.meta,positioning,calc,fit}
\usepackage{microtype}
\usepackage{accsupp,stringenc}
\newcommand{\copyformula}[2]{\begingroup
  \StringEncodingConvert\FormulaText{\detokenize{#1}}{utf8}{utf16be}%
  \BeginAccSupp{method=hex,unicode,ActualText={\pdfescapehex{\FormulaText}}}%
  #2\EndAccSupp{}\endgroup}
\ifdefined\pdfglyphtounicode
  \pdfglyphtounicode{Gamma}{0393}
  \pdfglyphtounicode{Sigma}{03A3}
  \pdfglyphtounicode{Pi}{03A0}
  \pdfglyphtounicode{rho}{03C1}
  \pdfglyphtounicode{pi}{03C0}
  \pdfglyphtounicode{sigma}{03C3}
  \pdfglyphtounicode{tfm:lmsy10/I}{2110}\pdfglyphtounicode{tfm:lmsy8/I}{2110}
  \pdfglyphtounicode{tfm:lmsy7/I}{2110}\pdfglyphtounicode{tfm:lmsy6/I}{2110}
  \newcommand{\mapmathsizes}[2]{%
    \pdfglyphtounicode{#1big}{#2}\pdfglyphtounicode{#1Big}{#2}%
    \pdfglyphtounicode{#1bigg}{#2}\pdfglyphtounicode{#1Bigg}{#2}}
  \mapmathsizes{parenleft}{0028}\mapmathsizes{parenright}{0029}
  \mapmathsizes{bracketleft}{005B}\mapmathsizes{bracketright}{005D}
  \mapmathsizes{braceleft}{007B}\mapmathsizes{braceright}{007D}
  \mapmathsizes{angbracketleft}{27E8}\mapmathsizes{angbracketright}{27E9}
  \mapmathsizes{slash}{002F}\mapmathsizes{radical}{221A}\pdfglyphtounicode{vextendsingle}{007C}
  \pdfglyphtounicode{summationtext}{2211}\pdfglyphtounicode{summationdisplay}{2211}
  \pdfglyphtounicode{producttext}{220F}\pdfglyphtounicode{productdisplay}{220F}
  \pdfglyphtounicode{circleplusdisplay}{2A01}
  \pdfglyphtounicode{hatwide}{0302}\pdfglyphtounicode{tildewide}{0303}
  \pdfglyphtounicode{tildewider}{0303}
\fi
\usepackage{fancyhdr}
\usepackage[section]{placeins}
\usepackage[hidelinks,pdfencoding=auto,psdextra]{hyperref}
\hypersetup{pdftitle={Minor--Skein Spatial Negami Splitting Vectors and a Gluing Formula for Marked State Sums},
  pdfauthor={Iwao Mizukai},
  pdfsubject={Spatial graphs, boundary states and gluing},
  pdfkeywords={spatial graph, Negami polynomial, graph minor, ribbon category, quantum group}}
\DeclareRobustCommand{\qMSNtag}{\ensuremath{\mathsf{qMSN}}}
\DeclareRobustCommand{\qMSNvector}[1]{\ensuremath{\mathbf{v}_{#1}^{\qMSNtag}}}
\newcommand{\edgepair}[2]{\mathopen{(}#1,\,#2\mathclose{)}}
\newcommand{\reserveblock}[1]{\par\ifdim\dimexpr\pagegoal-\pagetotal\relax<#1\relax\newpage\fi}
\newtheorem{theorem}{Theorem}[section]
\newtheorem{lemma}[theorem]{Lemma}
\newtheorem{proposition}[theorem]{Proposition}
\newtheorem{corollary}[theorem]{Corollary}
\theoremstyle{definition}
\newtheorem{definition}[theorem]{Definition}
\newtheorem{assumption}[theorem]{Assumption}

\theoremstyle{remark}
\newtheorem{remark}[theorem]{Remark}
\setlist{topsep=4pt,itemsep=2pt,parsep=0pt}
\newcommand{\papertitle}{Minor--Skein Spatial Negami Splitting Vectors\\[2pt] and a Gluing Formula for Marked State Sums}
\title{\LARGE\bfseries\papertitle}
\author{Iwao Mizukai\\[3pt] \normalsize Chiba Keizai University High School}
\date{\small Revision 9, 20 November 2026}
\begin{document}
\maketitle
\begin{abstract}
We organize the edge states of a decorated spatial graph by the minors
of their underlying abstract graphs. For pieces with disjoint edge sets,
Negami's component formula and functorial ribbon evaluation express
each marked sector as a boundary pairing, provided the local data
determine its marker. We give two explicit rules for this last step.
A two-terminal rule for $3$-connected targets tests each piece with
and without a virtual edge. At arbitrary finite vertex boundaries,
records of partial branch sets determine the occurrence of each fixed
simple target. The rules assign a weight once per edge state, even
when the state contains several minor models.

Inserting a fixed edge in $K_7-e$ creates a Petersen-family certificate
in $2{,}220$ labelled states. For a specified pair of triangles in $K_6$,
the same minor sector distinguishes an unlink from a Hopf link through
its bracket coefficient. The scalar marking in both examples preserves
a product pairing. A separate, prescribed two-strand cable gives an
example with nonproportional quantum blocks. At four fundamental ports
and generic parameter, reversible transport on finite labels forces the
blocks to be proportional under full braid-slide covariance. For
undecorated graphs we also construct a finite coloring model using
spanning pairs of cycles. A minor bit selects a Fox coloring count,
and two graph states realize closure forms of ranks two and one.

\end{abstract}
\begin{quotation}\small
\end{quotation}

\section{Introduction}

Cut a decorated spatial graph into two pieces. To compute its state sum
from the pieces, one must know which components meet across the cut,
how the two ribbon diagrams compose, and which graph minors can be
completed by gluing. Negami's splitting theory organizes the first
calculation by partitions of the terminal vertices~\cite{negami87}.
The functor on colored ribbon graphs supplies the second
calculation~\cite{rt,turaev}. We use ordinary graph minors to answer
the third question for each edge state. The spatial data in the
evaluation include the embedding, colors, framing, and vertex coupons.

The missing edge $uv$ in $K_6-uv$ illustrates the problem. This graph
and the single edge $uv$ each have zero $K_6$ marker, but their union
has marker one. The contribution here is a description of boundary
data and updates that place each weighted state in the correct minor
sector. The occurrence test is Boolean even when several models give
the same certificate. Partial-model methods for graph minors have an
established algorithmic literature~\cite{adler}. We detail the connection
with those methods in Section~\ref{sec:four-terminal-models}.

\reserveblock{95pt}
\paragraph{Guide to the results.}
Sections~3--5 give the two inputs to the pairing. The partition join counts components
after gluing, and composition evaluates the glued ribbon diagram.
The additional decision in the state-sum formula of Theorem~\ref{thm:q-gluing} is which marker
variable multiplies a given summand. For simple $3$-connected targets on
at least four vertices, Proposition~\ref{prop:torso-marker} answers this
at two terminals. Proposition~\ref{prop:partial-branch-records} handles
fixed simple targets at any finite vertex boundary. For calculations,
Section~\ref{sec:spatial-constituent} evaluates a specified constituent
link, while Section~\ref{sec:complete-marker-census} counts the abstract
minor certificates of labelled $K_7$ states. Section~\ref{sec:cable-probe}
constructs a nonproduct pairing with a prescribed cable.
Proposition~\ref{prop:finite-return} gives an obstruction for generic
four-point Temperley--Lieb transport. Proposition~\ref{prop:fox-nonproduct}
constructs a nonproduct pairing on the graph's own cycles in a finite
coloring theory. The analogous construction in the generic ribbon
theory remains open.

\subsection{Compression into boundary states}

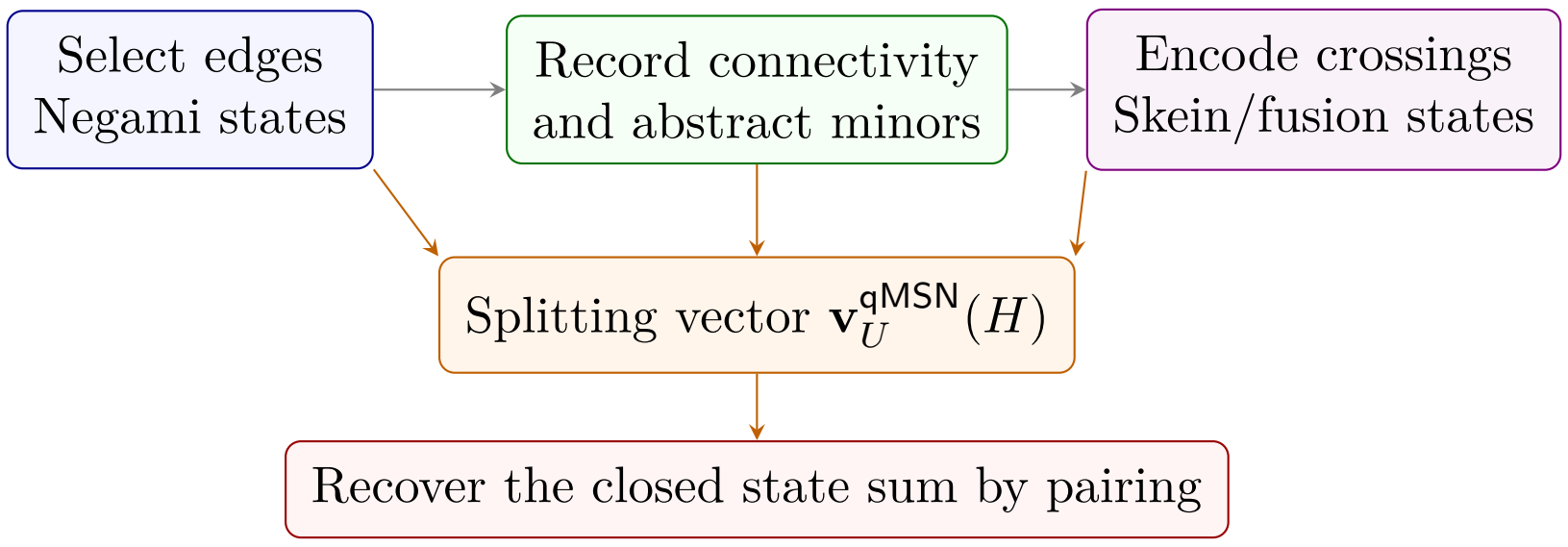
\begin{figure}[htbp]
\centering
\begin{tikzpicture}[x=1cm,y=1cm,
 panel/.style={draw,rounded corners=3pt,align=center,inner xsep=5pt,inner ysep=5pt,font=\small},
 stage arrow/.style={->,>=stealth,draw=black!50},
 merge/.style={->,>=stealth,draw=orange!75!black}]
\node[panel,draw=blue!55!black,fill=blue!4] (a) at (-3.9,0)
 {Select edges\\Negami states};
\node[panel,draw=green!45!black,fill=green!4] (b) at (0,0)
 {Record connectivity\\and abstract minors};
\node[panel,draw=violet,fill=violet!5] (c) at (3.9,0)
 {Encode crossings\\Skein/fusion states};
\node[panel,draw=orange!75!black,fill=orange!8] (v) at (0,-1.55)
 {Splitting vector $\qMSNvector{U}(H)$};
\node[panel,draw=red!60!black,fill=red!4] (z) at (0,-2.75)
 {Recover the closed state sum by pairing};
\draw[stage arrow] (a) -- (b);
\draw[stage arrow] (b) -- (c);
\draw[merge] (a.south east) -- (v.north west);
\draw[merge] (b) -- (v);
\draw[merge] (c.south west) -- (v.north east);
\draw[merge] (v) -- (z);
\end{tikzpicture}
\caption{Each edge state supplies an abstract boundary record and a ribbon evaluation. Pairing the vectors recovers the
closed state sum.}\label{fig:qmsn-compression}
\end{figure}

Each edge of a piece $H$ can be retained or deleted. There are therefore
$2^{|E(H)|}$ edge states:
\[
    \bigl|\{A\subseteq E(H)\}\bigr| = 2^{\lvert E(H)\rvert}.
\]
After cutting along a finite boundary $U$, we retain in the boundary
vector of each piece just the information needed for gluing.
Figure~\ref{fig:qmsn-compression} illustrates the construction.

\begin{remark}\label{stmt:1.1}
Splitting is realized by a matrix pairing of boundary state vectors.
For graph families with a fixed boundary type, this gives a
finite-dimensional transfer-matrix computation.
\end{remark}

\subsection{Notation}

The superscript \textsf{qMSN} means quantum Minor--Skein Negami.
Pieces $H,K$ share the boundary $U$ and have edge states $A,B$, respectively.

\begin{center}
\begingroup
\renewcommand{\arraystretch}{1.12}
\begin{tabularx}{\linewidth}{@{}l>{\raggedright\arraybackslash}X@{}}
\toprule
\textbf{Local notation} & \textbf{Meaning}\\
\midrule
$H_A,\; k_U^0(A)$ & Spanning subgraph on retained edges $A$, and its count of components disjoint from $U$\\
$\pi_A,\; \pi_A\vee\pi_B$ & Boundary connectivity partition, and its join with the partition from $K_B$\\
$t,\;X,\;Y$ & Component, retained-edge, and deleted-edge variables. The local weight is $t^{k_U^0(A)}X^{|A|}Y^{|E(H)|-|A|}$.\\
$F_q,\;F_q(H_A)$ & Ribbon-tangle functor and its value on the decorated state (Definition~\ref{stmt:5.3})\\
$q,\; A,\;s$ & Quantum-group parameter in $U_q(\mathfrak{sl}_2)$, bracket variable, and Jones variable $s=A^{-4}$, respectively. The bracket variable $A$ is a scalar.\\
$\qMSNvector{U}(H)$ & Weighted boundary vector (Definition~\ref{def:qmsn-vector})\\
\bottomrule
\end{tabularx}
\endgroup
\end{center}

\begin{center}
\begingroup
\renewcommand{\arraystretch}{1.12}
\begin{tabularx}{\linewidth}{@{}l>{\raggedright\arraybackslash}X@{}}
\toprule
\textbf{Closed notation} & \textbf{Meaning}\\
\midrule
$m(C)=(i,j)$ & IL and IK certificate bits for the abstract graph of a closed state $C$, relative to the chosen catalogues (Definition~\ref{stmt:2.7})\\
$\boldsymbol z=(z_{00},z_{01},z_{10},z_{11})$ & Four independent formal variables. A state with $m(C)=(i,j)$ receives $z_{ij}$.\\
$Z^{\qMSNtag}(G;\boldsymbol z)$ & Marked sum of the admissible states, with the chosen ribbon evaluations (Definition~\ref{def:closed-qmsn})\\
$[z_{ij}]Z^{\qMSNtag}(G;\boldsymbol z)$ & Sum of the contributions in the single sector $(i,j)$\\
$Z^{\qMSNtag}(G)$ & Unmarked sum obtained by setting all four $z_{ij}=1$, retaining the ribbon evaluations\\
$Z_{\mathrm{mark}}(G)$ & Shorthand in the census of Section~\ref{sec:complete-marker-census} for the marked sum with tensor-unit decoration, the full Petersen family as IL targets, and $\{K_7\}$ as the IK catalogue\\
\bottomrule
\end{tabularx}
\endgroup
\end{center}

\section{Spatial graphs and minor markers}

\begin{definition}[Spatial graph]\label{stmt:2.1}
A \emph{spatial graph} is the image
\[
  \Gamma=f(G)\subset S^3
\]
of a piecewise-linear embedding $f\colon G\hookrightarrow S^3$ of a
finite graph $G$. Thus $G$ denotes the abstract graph and $\Gamma$
its realization in three-dimensional space.
\end{definition}

\begin{definition}[Constituent knots and links]\label{stmt:2.2}
For a cycle $\gamma$ of $G$, the knot $f(\gamma)$ is a
\emph{constituent knot}. For a pair of vertex-disjoint cycles
$\lambda=\gamma_1\sqcup\gamma_2$, the link $f(\lambda)$ is a
\emph{constituent two-component link}.
\end{definition}

\begin{definition}[Intrinsic linking and knotting]\label{stmt:2.3}
A graph $G$ is \emph{intrinsically linked} (IL) if every spatial
embedding of $G$ contains a nonsplit link. It is \emph{intrinsically
knotted} (IK) if every spatial embedding of $G$ contains a nontrivial
knot.
\end{definition}

\begin{remark}\label{stmt:2.4}
Detecting a knot in a particular embedding establishes an existence
statement. Intrinsic knotting quantifies over all embeddings. Thus a
Jones-polynomial calculation, for example, that detects a knot in a
single embedding cannot by itself establish intrinsic knotting.
\end{remark}

\begin{definition}[Underlying abstract graph and graph minor]\label{stmt:2.5}
Given $\Gamma=f(G)$ and $C\subseteq E(G)$, write $\Gamma_C=f(G_C)$ and put
\[
  \operatorname{Abs}(\Gamma_C)=G_C=(V(G),C).
\]
Here $\operatorname{Abs}$ forgets the embedding and ribbon decorations.
It retains the specified vertices, including isolated vertices and labelled
boundary vertices. For a graph already presented abstractly, $\operatorname{Abs}$
acts as the identity. A graph $H$ is a \emph{minor} of an abstract graph
$G$, written $H\preceq G$, if it can
be obtained from $G$ by a finite sequence of vertex deletions, edge
deletions, and edge contractions. We contract edges with distinct
endpoints. A spatial-minor question would prescribe an embedded target
and allowable moves on embeddings. Our relation $\preceq$ compares
abstract graphs up to isomorphism.
\end{definition}

\begin{proposition}[Upward monotonicity of IL and IK certificates]\label{stmt:2.6}
Suppose that $H\preceq G$.
\begin{enumerate}
\item If $H$ is IL, then $G$ is IL.
\item If $H$ is IK, then $G$ is IK.
\end{enumerate}
\end{proposition}

\begin{proof}
Suppose that $G$ admits a linkless embedding or a knotless embedding,
respectively. Deletions can be performed directly in the embedding.
An edge contraction can be realized locally by enclosing the edge in a
small regular neighborhood that is a ball and shrinking that ball to a
point. The resulting minor therefore has a linkless or knotless embedding,
respectively. Applying this observation to $H$ gives the two assertions by
contraposition.
\end{proof}

\begin{definition}[Minor marker]\label{stmt:2.7}
Fix a set $\mathcal P_{\mathrm{IL}}$ of obstructions to linkless embedding
and a set $\mathcal P_{\mathrm{IK}}^{\mathrm{known}}$ of known obstructions
to knotless embedding. For an edge state $A$, define
\[
  m(A)=\bigl(m_{\mathrm{IL}}(A),m_{\mathrm{IK}}(A)\bigr)\in\{0,1\}^2,
\]
where
\begin{align*}
m_{\mathrm{IL}}(A)
  &=\mathbf 1[\exists P\in\mathcal P_{\mathrm{IL}}\colon P\preceq\operatorname{Abs}(H_A)],\\
m_{\mathrm{IK}}(A)
  &=\mathbf 1[\exists Q\in\mathcal P_{\mathrm{IK}}^{\mathrm{known}}
                       \colon Q\preceq\operatorname{Abs}(H_A)].
\end{align*}
The target catalogues are fixed before summing. For a closed state $C$ use the same tests on $G_C=\operatorname{Abs}(f(G_C))$ to define $m(C)$. The target has named vertices for the branch-set test and carries no prescribed knot or link type. Proposition~\ref{prop:torso-marker} computes these closed bits from refined two-terminal data for $3$-connected targets. Proposition~\ref{prop:partial-branch-records} gives the general finite-boundary rule for simple targets. The marker records occurrence, so a state containing several models of the same target is counted once. Its weight contributes to $[z_{ij}]Z^{\qMSNtag}(G;\boldsymbol z)$ for $(i,j)=m(C)$ (Definition~\ref{def:closed-qmsn}). When a witness is required, the branch sets and their connecting edges may also be retained.
\end{definition}

\begin{remark}\label{stmt:2.8}
With all seven Petersen-family targets, $m_{\mathrm{IL}}=1$ if and
only if the underlying graph is intrinsically linked~\cite{rst}. With a
smaller IL catalogue, the bit tests only its listed targets. Likewise
$m_{\mathrm{IK}}=1$ certifies intrinsic knotting, while $m_{\mathrm{IK}}=0$
means that none of the specified IK targets occurs. Knotless embeddability
is minor-closed by Proposition~\ref{stmt:2.6}, so the graph-minor
theorem~\cite{gmxx} gives a finite obstruction set. Our computations
use explicitly named targets.

For example, take the single IK target $K_7$. Foisy proved that
$K_{3,3,1,1}$ is intrinsically knotted~\cite{foisy}, yet its $K_7$ bit
is zero. Seven branch sets on eight vertices would include at least six
singletons, pairwise adjacent in a $K_7$ model. A clique in $K_{3,3,1,1}$
has at most one vertex from each of its four parts. For comparison, if
the two triples are $a_1,a_2,a_3$ and $b_1,b_2,b_3$, contracting
$a_1b_1$ and $a_2b_2$ produces a $K_6$ after parallel edges are
suppressed. Thus even the marker $(1,0)$ occurs on an intrinsically
knotted graph for the respective singleton catalogues $\{K_6\},\{K_7\}$.
\end{remark}

\subsection{Connection with the spatial Ramsey theorem}

\begin{definition}[Rectilinear spatial embedding]\label{stmt:2.9}
A spatial embedding is \emph{linear} or \emph{rectilinear} if it maps
each edge to a straight line segment in $\mathbb R^3$.
\end{definition}

\begin{theorem}[Negami's spatial Ramsey theorem]\label{stmt:2.10}
Let $G_1,\ldots,G_s$ be a finite list of spatial graphs. There is an
integer $R(G_1,\ldots,G_s)$ with the following property.
For $n>R(G_1,\ldots,G_s)$, given any rectilinear embedding of $K_n$ and
any edge coloring
\[
   E(K_n)=E_1\sqcup\cdots\sqcup E_s,
\]
there are an index $i$ and a subgraph $H$ such that $E(H)\subseteq E_i$
and the embedded subgraph $H$ is ambient isotopic to a subdivision of $G_i$.
\end{theorem}

\begin{proof}[Outline of Negami's proof]
The proof has two steps. First, subdivide the edges of each $G_i$ by
adding vertices of degree two as needed, and construct a standard diagram
$G_i^+$ with its vertices arranged in a line. Two standard diagrams
whose endpoint orders and projected crossing signs agree represent
ambient isotopic spatial graphs.

Second, start with a sufficiently large point set in general position
and apply the finite Ramsey theorem repeatedly. This gives a sequence
of points on which the endpoint order, the crossing data determined by
each quadruple, and the prescribed edge colors are uniform on the
required subsets. The resulting monochromatic complete subgraph contains
a subdivision with the endpoint order and crossing signs of $G_i^+$ for
the corresponding $i$. By the first step, it is ambient isotopic to a
subdivision of $G_i$. Only finitely many Ramsey extractions are needed,
so the required bound $R(G_1,\ldots,G_s)$ is finite.
See Negami~\cite{negami91} for the combinatorial bounds and the
construction of the standard diagrams.
\end{proof}

\begin{remark}\label{stmt:2.11}
The rectilinearity hypothesis is essential to this application of the
theorem. Here the Ramsey theorem guarantees the appearance of a prescribed
type in a sufficiently large rectilinearly embedded complete graph.
The marked sum below counts edge states with a specified abstract certificate. The Ramsey statement concerns the embedded type of a subgraph and supplies context for that distinction.
\end{remark}

\section{Boundaries, edge states, and set partitions}

\begin{definition}[Spatial graph tangle with boundary $U$]\label{stmt:3.1}
Let $f$ embed a finite graph $H$ in a $3$-manifold $M$, and let
$\copyformula{Γ}{\Gamma}=f(H)$ denote its image, as in Definition~\ref{stmt:2.1}.
Here $M$ is the manifold of the piece $H$, also denoted by $M_H$.
Let $\copyformula{Σ}{\Sigma}$ be the oriented cutting surface, viewed on this side of the
cut as a subsurface of $\partial M$, and let $U=\{u_1,\ldots,u_r\}$ be a
finite set of marked points on $\copyformula{Σ}{\Sigma}$. The embedded graph is a
\emph{spatial graph tangle with boundary $U$} if
\[
   \copyformula{Γ}{\Gamma}\cap\partial M = \copyformula{Γ}{\Gamma}\cap\copyformula{Σ}{\Sigma} = U
\]
and each edge incident with $U$ meets $\copyformula{Σ}{\Sigma}$ transversely there.
Thus all intersections with the manifold boundary occur at the marked
points of the cutting surface. When the embedding is understood, we
write $H$ for its image. Gluing a second piece $K$ along $(\copyformula{Σ}{\Sigma},U)$
then gives the embedded graph $H\cup_U K$.
\end{definition}

\begin{definition}[Edge state and spanning convention]\label{stmt:3.2}
An \emph{edge state} is a subset $A\subseteq E(H)$. The spanning
subgraph $H_A$ retains all vertices of $H$ and has edge set $A$.
A vertex incident with no selected edge counts as an isolated vertex.
\end{definition}

\begin{definition}[Interior component count and boundary partition]\label{stmt:3.3}
Let $k_U^0(A)$ be the number of connected components of $H_A$ disjoint
from $U$. Let $\pi_A\in\Pi(U)$ be the partition of $U$ in which $u_i$
and $u_j$ belong to the same block exactly when they belong to the same
connected component of $H_A$.
\end{definition}

\begin{definition}[Join of partitions]\label{stmt:3.4}
Let $\copyformula{Π}{\Pi}(U)$ denote the set of partitions of $U$. For $\copyformula{ρ}{\rho}\in\copyformula{Π}{\Pi}(U)$,
write $R(\copyformula{ρ}{\rho})$ for the equivalence relation on $U$ in which two points
are equivalent when they lie in the same block of $\copyformula{ρ}{\rho}$.
The \emph{join} of $\copyformula{π}{\pi},\copyformula{σ}{\sigma}\in\copyformula{Π}{\Pi}(U)$ is the partition $\copyformula{π}{\pi}\vee\copyformula{σ}{\sigma}$
whose equivalence relation is
\[
   R(\copyformula{π}{\pi}\vee\copyformula{σ}{\sigma}) = \operatorname{TC}(R(\copyformula{π}{\pi})\cup R(\copyformula{σ}{\sigma})),
\]
where $\operatorname{TC}$ denotes transitive closure. Equivalently, two
points lie in the same block of the join exactly when a finite chain
connects them, each step staying in a block of $\copyformula{π}{\pi}$ or in a block of
$\copyformula{σ}{\sigma}$. This gives the finest common coarsening of $\copyformula{π}{\pi}$ and $\copyformula{σ}{\sigma}$.
The number of blocks is denoted by $\lvert\copyformula{π}{\pi}\vee\copyformula{σ}{\sigma}\rvert$.
\end{definition}

\begin{figure}[htbp]
\centering
\begin{tikzpicture}[x=1cm,y=.76cm,font=\small,>=stealth]
\filldraw[rounded corners=3pt,draw=blue!55!black,fill=blue!4]
  (-4.3,0.45) rectangle (-1.48,-2.95);
\filldraw[rounded corners=3pt,draw=green!42!black,fill=green!4]
  (1.48,0.45) rectangle (4.3,-2.95);
\node at (-2.9,0.82) {$\pi=\{1,2\}\mid\{3,4\}$};
\node at (2.9,0.82) {$\sigma=\{1,3\}\mid\{2,4\}$};
\foreach \s/\t in {0/1,2/3}
 {\draw[thick,blue!55!black] (-1.6,{-.25-.76*\s})
  .. controls (-3.22,{-.25-.76*\s-.08}) and (-3.22,{-.25-.76*\t+.08})
  .. (-1.6,{-.25-.76*\t});}
\foreach \s/\t in {0/2,1/3}
 {\draw[thick,green!42!black] (1.6,{-.25-.76*\s})
  .. controls (3.2,{-.25-.76*\s-.08}) and (3.2,{-.25-.76*\t+.08})
  .. (1.6,{-.25-.76*\t});}
\foreach \j in {1,2,3,4}{
 \pgfmathsetmacro{\vy}{-.25-.76*(\j-1)}
 \filldraw[draw=orange!80!black,fill=white] (-1.6,\vy) circle (1.9pt);
 \node[right] at (-1.4,\vy) {$u_{\j}$};
 \filldraw[draw=orange!80!black,fill=white] (1.6,\vy) circle (1.9pt);
 \node[left] at (1.4,\vy) {$u_{\j}$};}
\node[draw=orange!75!black,fill=orange!8,rounded corners=3pt,
 align=center,inner sep=4pt] (join) at (0,-4.3)
 {$\pi\vee\sigma=\{1,2,3,4\}$\\Alternating between the two sides connects all four points.};
\draw[->,orange!75!black] (-1.5,-2.98) -- ([xshift=-5pt]join.north);
\draw[->,orange!75!black] (1.5,-2.98) -- ([xshift=5pt]join.north);
\end{tikzpicture}
\caption{The labels $1,2,3,4$ stand for $u_1,u_2,u_3,u_4$. The join records connectivity after imposing the connections on both sides.}\label{fig:qmsn-join}
\end{figure}
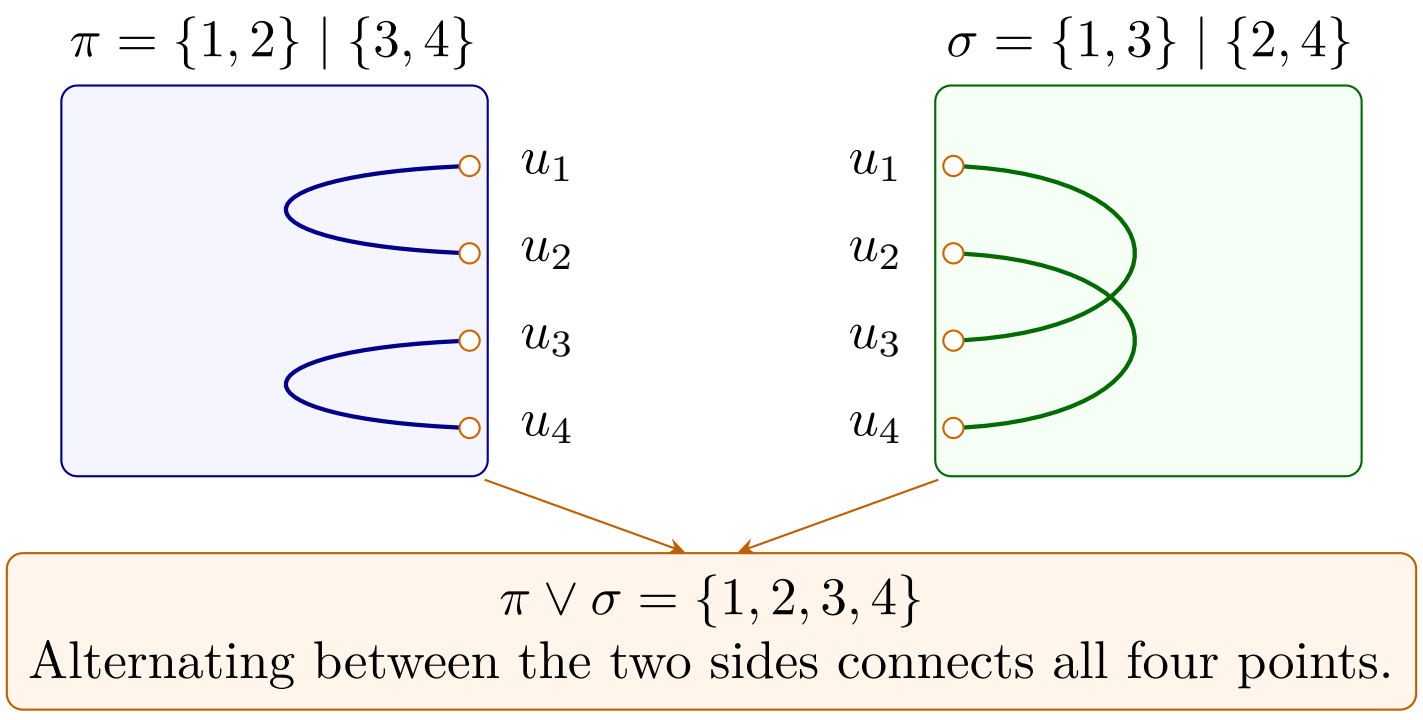

\begin{lemma}[Join formula for the number of components]\label{stmt:3.5}
For $A\subseteq E(H)$ and $B\subseteq E(K)$,
\[
  k(A\cup_U B)=k_U^0(A)+k_U^0(B)+|\pi_A\vee\pi_B|.
\]
\end{lemma}

\begin{proof}
Each component of $H_A$ disjoint from $U$ remains a separate component
after gluing, since it cannot meet $K_B$. There are $k_U^0(A)$ such
components on the $H$ side and $k_U^0(B)$ on the $K$ side.
For the remaining components, two boundary points become connected
if they can be joined within $H_A$, within $K_B$, or by a finite
sequence of paths alternating between the two. This is precisely the
equivalence relation of $\pi_A\vee\pi_B$. Its $|\pi_A\vee\pi_B|$
blocks count the components meeting the boundary after gluing. The
three classes of components are disjoint, giving the formula.
\end{proof}

\subsection{The fifteen partitions of a four-point boundary}

Place $U=\{1,2,3,4\}$ on a circle in this order. The fourth Bell
number is $B_4=15$. A partition is \emph{noncrossing} if, when the points of
each block are joined by chords inside the disk, chords belonging to
different blocks do not cross.
Each row gives one partition. Braces enclose a block, and a vertical
bar separates adjacent blocks. The count refers to all the rows of the
same block sizes and crossing type.
\begin{center}
\begingroup
\renewcommand{\arraystretch}{1.10}
\setlength{\tabcolsep}{7pt}
\begin{tabular*}{\linewidth}{@{\extracolsep{\fill}}lcll@{}}
\toprule
\textbf{Block sizes} & \textbf{Count} & \textbf{Partition of $\{1,2,3,4\}$} & \textbf{Type}\\
\midrule
$1+1+1+1$ & 1 & $\{1\}\;\vert\;\{2\}\;\vert\;\{3\}\;\vert\;\{4\}$ & Noncrossing\\
\addlinespace[2pt]
$2+1+1$ & 6 & $\{1,2\}\;\vert\;\{3\}\;\vert\;\{4\}$ & Noncrossing\\
 & & $\{1,3\}\;\vert\;\{2\}\;\vert\;\{4\}$ & \\
 & & $\{1,4\}\;\vert\;\{2\}\;\vert\;\{3\}$ & \\
 & & $\{2,3\}\;\vert\;\{1\}\;\vert\;\{4\}$ & \\
 & & $\{2,4\}\;\vert\;\{1\}\;\vert\;\{3\}$ & \\
 & & $\{3,4\}\;\vert\;\{1\}\;\vert\;\{2\}$ & \\
\addlinespace[2pt]
$2+2$ & 2 & $\{1,2\}\;\vert\;\{3,4\}$ & Noncrossing\\
 & & $\{1,4\}\;\vert\;\{2,3\}$ & \\
\addlinespace[2pt]
$2+2$ & 1 & $\{1,3\}\;\vert\;\{2,4\}$ & \textbf{Crossing}\\
\addlinespace[2pt]
$3+1$ & 4 & $\{1,2,3\}\;\vert\;\{4\}$ & Noncrossing\\
 & & $\{1,2,4\}\;\vert\;\{3\}$ & \\
 & & $\{1,3,4\}\;\vert\;\{2\}$ & \\
 & & $\{2,3,4\}\;\vert\;\{1\}$ & \\
\addlinespace[2pt]
$4$ & 1 & $\{1,2,3,4\}$ & Noncrossing\\
\bottomrule
\end{tabular*}
\endgroup
\end{center}

\begin{proposition}[Crossing partitions of a four-point boundary]\label{stmt:3.6}
For the cyclic order $1,2,3,4$, the only crossing partition among the
fifteen is $\{1,3\}\mid\{2,4\}$. The other fourteen are noncrossing.
\end{proposition}

\begin{proof}
A crossing requires two points $a,c$ in one block and two points $b,d$
in another, alternating in cyclic order as $a,b,c,d$. With only four
points this forces the partition $\{1,3\}\mid\{2,4\}$. Conversely, the
chord joining $1$ to $3$ crosses the chord joining $2$ to $4$ in the disk. Each remaining partition in
the table admits a drawing without crossings between blocks, using
chords or a single polygon.
\end{proof}

\section{The classical Negami boundary space}

Take the coefficient ring to be
\[
 R_0=\mathbb Z[t,X,Y].
\]
When quantum parameters enter, we extend scalars to a commutative ring
$R$ containing $R_0$.

\begin{definition}[Partition module]\label{stmt:4.1}
The \emph{classical partition module} for the boundary $U$ is
\[
 \mathcal N_U\coloneqq\bigoplus_{\pi\in\Pi(U)}R_0e_\pi.
\]
\end{definition}

\begin{definition}[Negami boundary pairing]\label{stmt:4.2}
On basis vectors, set
\[
 \langle e_\pi,e_\sigma\rangle_U^N\coloneqq t^{|\pi\vee\sigma|},
\]
and extend bilinearly.
\end{definition}

\begin{definition}[Classical splitting vector]\label{stmt:4.3}
The classical splitting vector is
\[
 \mathbf v_U^N(H)\coloneqq
  \sum_{A\subseteq E(H)}t^{k_U^0(A)}X^{|A|}Y^{|E(H)|-|A|}
                             e_{\pi_A}\in\mathcal N_U.
\]
\end{definition}

\begin{remark}\label{stmt:4.4}
A component of $H_A$ disjoint from $U$ contributes one factor of $t$
to the local weight. If there are $a=\copyformula{k_U^0(A)}{k_U^0(A)}$ such
components, their combined contribution is $\copyformula{t^a}{t^a}$. The remaining components
meet the boundary, where some of them can merge through $K_B$. The basis
vector $e_{\pi_A}$ records their boundary partition. After gluing, the
components meeting $U$ contribute $\copyformula{t^b}{t^b}$, with one factor of $t$ for each
of the $b=\copyformula{|π_A ∨ π_B|}{\lvert\pi_A\vee\pi_B\rvert}$
blocks of the join.
\end{remark}

To compare with Negami's quotient coordinates~\cite[Theorem~4.2]{negami87},
write $a_\pi(H)$ for the coefficient of $e_\pi$ in $\mathbf v_U^N(H)$.
Write $Z^N$ for the ordinary Negami polynomial.
Let $H/\sigma$ identify the vertices in each block of $\sigma$, retaining
the edge set, including loops and parallel edges created by identification.
Then
\[
  Z^N(H/\sigma)=\sum_{\pi\in\Pi(U)} t^{|\pi\vee\sigma|}a_\pi(H).
\]
Indeed, identifying terminals changes the components meeting $U$
according to $\pi\vee\sigma$. In these direct state coordinates the
classical gluing matrix is $\bigl(t^{|\pi\vee\sigma|}\bigr)_{\pi,\sigma}$.
The marked construction refines each coordinate by its boundary data.

\section{Boundary states from quantum groups and ribbon categories}

\begin{definition}[Ribbon category data]\label{stmt:5.1}
Let $\Bbbk$ be a field and let $\mathcal C$ be a $\Bbbk$-linear ribbon category.
Assume $\operatorname{End}_{\mathcal C}(\mathbf 1)=\Bbbk$. Take $R$ to be a commutative $\Bbbk$-algebra compatible with the polynomial variables and extend the quantum morphism spaces to $R$ when forming tensor products in Section~6. Its tensor product, duals, braiding, evaluation and coevaluation maps, and twist
satisfy the ribbon category axioms. Assign an object $V_i$ and an orientation to
each boundary point $u_i$. At each interior vertex of valence $d$, choose a
morphism, represented by a coupon.
\end{definition}

\begin{definition}[Quantum boundary spaces]\label{stmt:5.2}
Write $V_U=V_1^{\epsilon_1}\otimes\cdots\otimes V_r^{\epsilon_r}$ for the tensor
product determined by the boundary orientations, and set
\[
  \mathcal Q_U^+=\operatorname{Hom}_{\mathcal C}(\mathbf 1,V_U),
  \qquad
  \mathcal Q_U^-=\operatorname{Hom}_{\mathcal C}(V_U,\mathbf 1).
\]
When the chosen boundary identifications allow both spaces to be written in the same coordinates, we write simply
$\mathcal Q_U$.
\end{definition}

\begin{definition}[Quantum evaluation functor]\label{stmt:5.3}
Let $F_q$ denote a monoidal functor from the category of ribbon tangles with
coupons to $\mathcal C$. It assigns $F_q(H_A)\in\mathcal Q_U^+$ to a state on
the left and $F_q(K_B)\in\mathcal Q_U^-$ to a state on the right. Choose the coupon for every allowed incidence pattern at a vertex and specify the quantum boundary object even at an unoccupied port. One may use a direct sum indexed by boundary occupations, with zero pairing between incompatible sectors. Section~\ref{subsec:explicit-contraction} instead gives explicit terminations at the tensor unit and hence a fixed two-strand boundary object.
\end{definition}

\begin{definition}[Quantum gluing evaluation]\label{stmt:5.4}
For $\alpha\in\mathcal Q_U^+$ and $\beta\in\mathcal Q_U^-$, put
\[
  \operatorname{ev}_U(\alpha,\beta)
   \coloneqq \beta\circ\alpha\in\operatorname{End}_{\mathcal C}(\mathbf 1)
   \cong\Bbbk .
\]
\end{definition}

\begin{assumption}[Local functoriality]\label{stmt:5.5}\label{ass:local-functoriality}
For each gluing of states,
\[
  F_q(H_A\cup_U K_B)
   =\operatorname{ev}_U\bigl(F_q(H_A),F_q(K_B)\bigr).
\]
For a ribbon functor of Reshetikhin--Turaev type, this is the basic property
that gluing diagrams corresponds to composing morphisms.
\end{assumption}

\begin{remark}\label{stmt:5.6}
The choice of a coupon is part of the data at a vertex of a spatial graph.
One must specify which morphism the vertex represents. Different choices of
coupons can give different invariants.
\end{remark}

\subsection{Fusion channels}

Consider the representation category of $U_q(\mathfrak{sl}_2)$ at generic
$q$. Denote the fundamental representation by the color $1$. Since
\[
  1\otimes 1\cong 0\oplus 2,
\]
the first fusion of two strands has two channels, distinguished by the
intermediate color $0$ or $2$. Six strands have five fusion paths ending in
total color $0$, namely
\[
  \begin{aligned}
    \Phi_{AA}&\colon 0,1,0,1,0,1,0, &
    \Phi_{AB}&\colon 0,1,0,1,2,1,0,\\
    \Phi_{BA}&\colon 0,1,2,1,0,1,0, &
    \Phi_{BB}&\colon 0,1,2,1,2,1,0,\\
    \Phi_{CC}&\colon 0,1,2,3,2,1,0. &&
  \end{aligned}
\]

\begin{figure}[htbp]
\centering
\begin{tikzpicture}[x=.95cm,y=.65cm,font=\small,line width=.8pt]
\draw (-1,0)--(6.2,0);
\node[left] at (-1,0) {$a_0=0$};
\node[right] at (6.2,0) {$a_6=0$};
\foreach \j in {0,...,5} {
  \draw[blue!65!black] (\j,0)--(\j,1);
  \fill[blue!65!black] (\j,0) circle (1.5pt);
  \node[above] at (\j,1) {$1$};
}
\foreach \j in {1,...,5} {\node[below] at (\j-.5,0) {$a_{\j}$};}
\end{tikzpicture}
\caption{Sequential fusion of six fundamental colors. Each upper edge carries color $1$. An admissible list of intermediate labels $a_i$ gives one of the five paths above.}

\label{fig:six-fusion}
\end{figure}
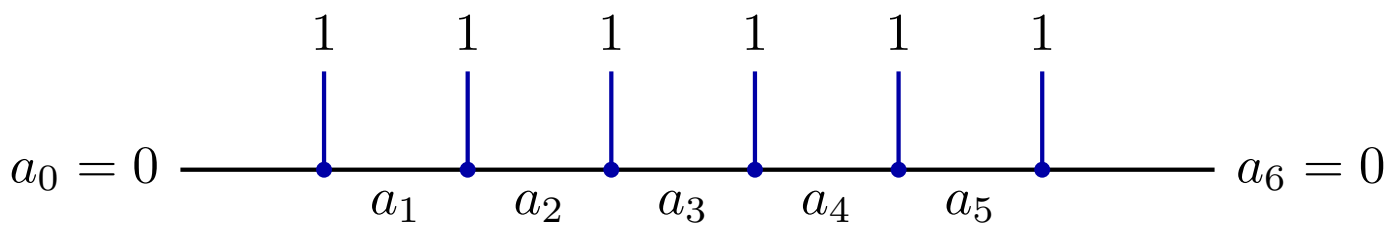

\begin{remark}\label{stmt:5.7}
For four points of fundamental color $1$ and total color $0$, the invariant
subspace has dimension two, with intermediate channels $0$ and $2$.
If the four points are grouped in pairs into two color-$2$ cables, then
$2\otimes 2$ decomposes into three channels of colors $0$, $2$, and $4$.
Whenever the three basis elements $P_0,P_2,P_4$ are used, the color-$2$
cabling convention must therefore be specified.
\end{remark}

\section{The Minor--Skein Spatial Negami state space}

\begin{definition}[Complete boundary datum]\label{stmt:6.1}
Let $\mathcal M_U$ be a finite set of boundary data. The datum $\mu_U(A)$
is computed from the piece $H_A$ alone. It includes its partition
$p(\mu)=\pi_\mu\in\Pi(U)$ and may include colors, orientations,
occupancies, and partial minor certificates. Fix a boundary test
$\operatorname{comp}(\mu,\nu)\in\{0,1\}$ for matching the two pieces.

A \emph{marker rule} for the chosen IL and IK catalogues is a function
$g_U$ from compatible pairs of data to $\{0,1\}^2$. It is \emph{exact}
if for every compatible pair of states in the class of gluings under
consideration it satisfies
\[
 g_U\bigl(\mu_U(A),\mu_U(B)\bigr)=m(A\cup_U B),
\]
where the right side is the graph-minor test of Definition~\ref{stmt:2.7}
on the glued state. This equality tests the proposed rule. Its inputs
are the two local data. For a fixed simple target $P$, Section~\ref{sec:four-terminal-models} includes the set $\Sigma_P(H_A;U)$ of partial branch-set records in $\mu_U(A)$. The rule tests whether at least one glued record completes a model of $P$.
\end{definition}

\begin{definition}[Classical marker space]\label{stmt:6.2}
Set
\[
 \widetilde{\mathcal N}_U\coloneqq\bigoplus_{\mu\in\mathcal M_U}Re_\mu,
 \qquad
 \langle e_\mu,e_\nu\rangle_U^{N}
   \coloneqq\operatorname{comp}(\mu,\nu)\,t^{|\pi_\mu\vee\pi_\nu|}.
\]
\end{definition}

\begin{definition}[$\qMSNtag$ boundary state spaces]\label{stmt:6.3}
Set
\[
 \mathcal S_U^+=\widetilde{\mathcal N}_U\otimes_R\mathcal Q_U^+,
 \qquad
 \mathcal S_U^-=\widetilde{\mathcal N}_U\otimes_R\mathcal Q_U^- .
\]
\end{definition}

\begin{definition}[$\qMSNtag$ splitting vector]\label{stmt:6.4}\label{def:qmsn-vector}
Let $\mu_U(A)$ be the complete boundary datum induced by the state $A$.
The splitting vector of the left-hand piece $H$ is
\[
 \qMSNvector{U}(H)
 =\sum_{A\subseteq E(H)}
   t^{k_U^0(A)}X^{|A|}Y^{|E(H)|-|A|}
    e_{\mu_U(A)}\otimes F_q(H_A).
\]
The right-hand splitting vector is defined similarly as an element of
$\mathcal S_U^-$.
\end{definition}

\begin{definition}[Product pairing]\label{stmt:6.5}
On pure tensors, set
\[
 \langle e_\mu\otimes\alpha,e_\nu\otimes\beta\rangle_U
 \coloneqq\operatorname{comp}(\mu,\nu)\,t^{|\pi_\mu\vee\pi_\nu|}
 \operatorname{ev}_U(\alpha,\beta)
\]
and extend bilinearly. For an exact rule $g_U$ define the marked pairing over $R[\boldsymbol z]$. On compatible basis data multiply the displayed value by $z_{g_U(\mu,\nu)}$, and set the marked value to zero on incompatible pairs. Denote this pairing by $\langle\ ,\ \rangle_{U,\boldsymbol z}$.
\end{definition}

\begin{definition}[Closed $\qMSNtag$ state sum]\label{stmt:6.6}\label{def:closed-qmsn}
For a closed decorated spatial graph $G$, let $\mathcal A(G)$ be its admissible edge subsets, specified by the chosen decoration. By default all edge subsets are admissible. With $w(C)=t^{k(C)}X^{|C|}Y^{|E(G)|-|C|}$ define
\[
 \begin{aligned}
 Z^{\qMSNtag}(G)&=\sum_{C\in\mathcal A(G)} w(C)F_q(G_C),\\
 Z^{\qMSNtag}(G;\boldsymbol z)&=\sum_{C\in\mathcal A(G)}
          w(C)\,z_{m(C)} F_q(G_C).
 \end{aligned}
\]
Here $z_{ij}$ are four independent indeterminates, with $z_{m(C)}=z_{ij}$ for $m(C)=(i,j)$. Thus $[z_{ij}]Z^{\qMSNtag}(G;\boldsymbol z)$ is the sum of states with that marker. Setting all $z_{ij}=1$ recovers the unmarked sum, with the edge, component and ribbon weights retained.

\end{definition}

\begin{remark}\label{stmt:6.7}
Each summand records the number of selected
edges, the number of components and their boundary connections, and the
quantum evaluation of the spatial crossings and fusions. The product
pairing combines these contributions by multiplication within each summand.
\end{remark}

\section{A gluing formula for marked state sums}
The identity in Theorem~\ref{thm:q-gluing} expresses a finite state sum
in boundary coordinates. Its classical component count comes from
Lemma~\ref{stmt:3.5}, its monomials multiply by Lemma~\ref{lem:q-weights},
and its ribbon factor composes by Lemma~\ref{lem:q-compose}. Given an
exact marker rule, distributivity groups these terms by their closed
marker. Sections~7.2 and~\ref{sec:four-terminal-models} construct such a
rule on two different choices of boundary data.
\subsection{State decomposition and the gluing identity}

\begin{lemma}[Unique decomposition of states]\label{stmt:7.1}\label{lem:q-decomp}
Suppose that $E(H)$ and $E(K)$ are disjoint. The map
\[
 \Psi\colon\mathcal P(E(H))\times\mathcal P(E(K))
  \longrightarrow\mathcal P(E(H\cup_U K)),\qquad
  (A,B)\longmapsto A\cup B,
\]
is a bijection.
\end{lemma}

\begin{proof}
Given a global state $C$, let $A=C\cap E(H)$ and $B=C\cap E(K)$.
Then $C=A\cup B$. Disjointness of the edge sets makes this decomposition
unique and hence gives the inverse map.
\end{proof}

\begin{lemma}[Multiplicativity of edge weights]\label{stmt:7.2}\label{lem:q-weights}
Let $G=H\cup_U K$ have disjoint piece edge sets. For a graph $L$ and
a state $S\subseteq E(L)$, write $\omega_L(S)$ for the product of one
factor $X$ for each retained edge and one factor $Y$ for each deleted edge. Thus
\[
 \copyformula{ω_L(S) = X^{|S|} Y^{|E(L)|−|S|}}
             {\omega_L(S)=X^{|S|}Y^{|E(L)|-|S|}} .
\]
For $A\subseteq E(H)$ and $B\subseteq E(K)$ the edge weight of their union is
\[
 \copyformula{ω_G(A ∪ B) = ω_H(A) ω_K(B)}
             {\omega_G(A\cup B)=\omega_H(A)\,\omega_K(B)} .
\]
\end{lemma}

\begin{proof}
Each retained edge lies in exactly one of $A$ and $B$. Likewise, each
deleted edge lies in exactly one of $E(H)\setminus A$ and $E(K)\setminus B$.
Thus the numbers of retained edges add, as do the numbers of deleted edges.
These are the exponents of $X$ and $Y$, respectively.
\end{proof}

\begin{lemma}[When the boundary data determine the marker]\label{stmt:7.3}
An exact marker rule exists on a class of gluings if and only if any two
compatible state pairs with the same ordered pair of boundary data have
the same global marker. For a general finite boundary, the sets of records in Proposition~\ref{prop:partial-branch-records} have this property for finite catalogues of simple targets. For $U=\{u,v\}$ and finite catalogues of simple
$3$-connected targets on at least four vertices, the following local
data have this property. For each piece $L$, record whether $u$ and $v$
are connected and, for each target $P$, whether $P$ is a minor of $L$
and of $L$ with the edge $uv$ added.
\end{lemma}

\begin{proof}
An exact rule must give the same value whenever its two inputs agree. Conversely, if the
marker is constant on each fiber, its value there defines the rule on
every realized pair of data. Unrealized pairs do not enter the sum.
For the stated two-terminal signature,
Proposition~\ref{prop:torso-marker} below gives the explicit Boolean
update for each target and proves it by restricting a minor model at
the separator. Taking the disjunction of the outputs within each
catalogue gives the two bits.
\end{proof}

\begin{lemma}[Composition of quantum evaluations]\label{stmt:7.4}\label{lem:q-compose}
Under Assumption~\ref{ass:local-functoriality},
\[
 \operatorname{ev}_U\bigl(F_q(H_A),F_q(K_B)\bigr)
   = F_q(H_A\cup_U K_B).
\]
\end{lemma}

\begin{proof}
The monoidal functor $F_q$ sends the gluing of tangles to composition in
$\mathcal C$. Composing the left morphism $\mathbf 1\to V_U$ with the right
morphism $V_U\to\mathbf 1$ evaluates the closed diagram.
\end{proof}

\begin{theorem}[Edge-disjoint gluing formula for marked state sums]\label{stmt:7.5}\label{thm:q-gluing}
Let $G=H\cup_U K$ with $V(H)\cap V(K)=U$ and disjoint piece edge sets. Assume local
functoriality as in Assumption~\ref{ass:local-functoriality}. Suppose
the union map of Lemma~\ref{lem:q-decomp} identifies compatible pairs
with admissible global states. More generally, it suffices that every compatible pair
glues to an admissible state and that every omitted admissible state $C$ has
$w(C)F_q(G_C)=0$. Assume first that the data admit an exact marker rule $g_U$.
Write the splitting vectors in the basis of boundary data as
\[
 \copyformula{v_U^{qMSN}(H) = ∑_μ e_μ ⊗ h_μ}
             {\qMSNvector{U}(H)=\sum_{\mu} e_\mu\otimes h_\mu},
 \qquad
 \copyformula{v_U^{qMSN}(K) = ∑_ν e_ν ⊗ k_ν}
             {\qMSNvector{U}(K)=\sum_{\nu} e_\nu\otimes k_\nu}.
\]
For each pair of bits $(i,j)$, define the sector index set by
\[
 \copyformula{ℐ_{ij}}{\mathcal I_{ij}} =
 \bigl\{(\mu,\nu)\in\mathcal M_U\times\mathcal M_U\!:
       \operatorname{comp}(\mu,\nu)=1,\ g_U(\mu,\nu)=(i,j)\bigr\}.
\]
Its elements are precisely the compatible ordered pairs whose glued
marker has IL bit $i$ and IK bit $j$. For a polynomial $P$ in the four
marker variables, $\copyformula{[z_{ij}] P}{[z_{ij}]P}$ denotes the
coefficient of the monomial $z_{ij}$. This sector coefficient is
\[
 \copyformula{[z_{ij}] Z^{qMSN}(G; z)}{[z_{ij}]\,Z^{\qMSNtag}(G;\boldsymbol z)}
 =
 \copyformula{∑_{(μ,ν)∈ℐ_{ij}}}{\sum_{(\mu,\nu)\in\mathcal I_{ij}}}
 \copyformula{t^{|π_μ ∨ π_ν|}}{t^{\lvert\pi_\mu\vee\pi_\nu\rvert}}\,
 \copyformula{ev_U}{\operatorname{ev}_U}(h_\mu,k_\nu) .
\]
Summing the four sector contributions with their respective marker variables gives
\[
 \copyformula{Z^{qMSN}(H ∪_U K; z) = ⟨v_U^{qMSN}(H), v_U^{qMSN}(K)⟩_{U,z}}
 {Z^{\qMSNtag}(H\cup_U K;\boldsymbol z)
 =\bigl\langle\qMSNvector{U}(H),
               \qMSNvector{U}(K)\bigr\rangle_{U,\boldsymbol z}} .
\]
For the two-terminal
catalogues of Lemma~\ref{stmt:7.3}, the signature and update of
Proposition~\ref{prop:torso-marker} supply the rule with no additional
marker-recovery assumption. At any fixed finite vertex boundary,
Proposition~\ref{prop:partial-branch-records} supplies the rule for every
fixed finite catalogue of simple targets. A computation that retains
coarser data must verify that those data still determine the marker.

Specializing every marker variable to one gives the unmarked identity
\[
 \copyformula{Z^{qMSN}(G) = ⟨v_U^{qMSN}(H), v_U^{qMSN}(K)⟩_U}
 {Z^{\qMSNtag}(G)
 =\bigl\langle\qMSNvector{U}(H),
               \qMSNvector{U}(K)\bigr\rangle_U} .
\]
This last identity also holds under the gluing hypotheses alone.
\end{theorem}

\begin{proof}
Expand the pairing using Definition~\ref{def:qmsn-vector}. Fix a
compatible pair $(A,B)$, set $C=A\cup_U B$, and write $m_H=|E(H)|$,
$m_K=|E(K)|$. Its contribution before assigning a sector is
\[
 \begin{aligned}
 &t^{k_U^0(A)+k_U^0(B)+|\pi_A\vee\pi_B|}
 X^{|A|+|B|}Y^{m_H+m_K-|A|-|B|}\,
 \operatorname{ev}_U\bigl(F_q(H_A),F_q(K_B)\bigr)\\
 &\hspace{2em}=w(C)F_q(G_C).
 \end{aligned}
\]
The three identifications are the component formula of
Lemma~\ref{stmt:3.5}, the edge weight identity of
Lemma~\ref{lem:q-weights}, and the quantum composition formula of
Lemma~\ref{lem:q-compose}. Each admissible state represented by a compatible pair occurs once by
Lemma~\ref{lem:q-decomp}. The hypothesis makes each omitted admissible
term zero, so summation proves the unmarked identity.

In the marked sum, the state $C$ contributes to the coefficient of
$z_{ij}$ precisely when $m(C)=(i,j)$. By exactness, this occurs precisely
when $(\mu_U(A),\mu_U(B))$ belongs to $\copyformula{ℐ_{ij}}{\mathcal I_{ij}}$.
Grouping first by these two data sums the local contributions to
$h_\mu$ and $k_\nu$ by bilinearity, and gives the displayed formula.
A zero closed contribution remains zero after multiplication by any
sector variable.
\end{proof}

\begin{remark}\label{stmt:7.6}
Theorem~\ref{thm:q-gluing} combines the classical splitting calculation
with the functorial evaluation of the specified ribbon diagrams. Its
marker hypothesis is a condition on what the two pieces remember.
Computing $g_U$ requires a separate graph-theoretic argument. Glue
$K_6-uv$ first to the edgeless graph on $u,v$, then to the single
edge $uv$. Both pairs have local $K_6$ bits $(0,0)$, while the two
closures have bits $0$ and $1$, respectively.
The unmarked identity only uses the state decomposition
and the classical and quantum composition rules.
\end{remark}

\subsection{Propagating minor certificates across a two-vertex boundary}

A minor model may use both sides of a cut. For example, neither
$K_6-uv$ nor the single edge $uv$ contains a $K_6$ minor, although their
union is $K_6$. We give a gluing rule for a finite catalogue
$\mathcal P$ of simple $3$-connected graphs, each on at least four
vertices. All minor tests in this subsection concern the underlying
abstract graph.

Let $U=\{u,v\}$, with $u\ne v$. For a graph $L$ with these terminals,
let $L^+$ be $L$ with the edge $uv$ added if it is absent. The added
edge is used only in the minor test and carries no state weight.
The restriction of a $3$-connected minor to an edge-completed side is
the familiar two-separation principle, used for example in the proof of
Scott and Sokal~\cite[Proposition~7.16]{scottsokal}. For spanning edge
states the other side may fail to connect the terminals. The bit $c$
below records when it actually supplies the missing path.
For $P\in\mathcal P$, put
\[
 a_P(L)=\mathbf 1[P\preceq L],\qquad
 b_P(L)=\mathbf 1[P\preceq L^+],\qquad
 c(L)=\mathbf 1[u\text{ and }v\text{ are connected in }L].
\]
Here $a_P\le b_P$, and all three flags are monotone under the addition of
edges. The bit $c(L)$ is already determined by the boundary partition.

\begin{proposition}[Two-terminal update for minor certificates]
\label{prop:torso-marker}
Let $L$ and $M$ have disjoint edge sets and exactly the vertices $u,v$ in
common. Write $J=L\cup_U M$. With $\lor$ and $\land$ denoting Boolean
disjunction and conjunction, respectively, the update is
\begin{align*}
 a_P(J)&= a_P(L)\lor a_P(M)
       \lor\bigl(b_P(L)\land c(M)\bigr)
       \lor\bigl(b_P(M)\land c(L)\bigr),\\
 b_P(J)&= b_P(L)\lor b_P(M),\qquad
 c(J)=c(L)\lor c(M).
\end{align*}
The signature consisting of $c$ and all pairs $(a_P,b_P)$ is therefore
closed under repeated gluing along $U$.
\end{proposition}

\begin{proof}
A model of $P$ consists of disjoint connected branch sets, one for each
vertex of $P$, with an edge witnessing each required adjacency.
A model on either side survives in $J$. If $M$ contains a $u$--$v$
path, deleting the rest of $M$ and contracting the path to a single
edge realizes $L^+$ as a minor of $J$ (discard a duplicate edge if
necessary). Thus $b_P(L)\land c(M)$ also suffices, as does the
symmetric condition.

Conversely, take a model of $P$ in $J$. At most two branch sets meet
$U$. Deleting their labels from $P$ leaves a connected, nonempty graph
by $3$-connectivity. Every remaining branch set lies on just one side, and
two such sets on opposite sides cannot be adjacent. They therefore all
lie on the same side, say $L$. Restrict each terminal branch set to $L$.
Each component of the restriction contains one of its terminals. Indeed, a path
from a vertex of the branch set to the first terminal encountered stays
in $L$ up to that terminal. Thus only a set containing both terminals
can split, and in that case the edge $uv$ reconnects it. If the
terminals lie in different branch sets, an adjacency between those two
sets witnessed only in $M$ is supplied by the same edge. Every required
adjacency to a branch set avoiding $U$ already lies in $L$. This gives
a model in $L^+$.

If $c(M)=0$, every component of $M$ meets at most one terminal.
An excursion of a branch set into $M$ can then reenter $L$ only at the
same terminal. Nor can an edge of $M$ join two different terminal branch
sets, since together with their paths to their respective terminals in $M$,
it would give a $u$--$v$ path. Hence in this case the restriction is a
model in $L$ itself. This proves the first formula. For the second,
apply the same restriction argument to $J^+$, assigning the extra edge
to either side. The two completed sides are still $L^+$ and $M^+$,
and both are subgraphs of $J^+$. Finally, a simple $u$--$v$ path lies
entirely on one side, proving the formula for $c$. Each update computes
the flags of the union itself, so repeated updates give the same result
in either order.
\end{proof}

For each side, include $c$ and all $(a_P,b_P)$ in its local datum. Applying the proposition target by target and taking disjunction within each of the IL and IK catalogues constructs the exact rule $g_U$ of Definition~\ref{stmt:6.1}. This supplies the hypothesis of the sector formula in Theorem~\ref{thm:q-gluing}. More generally, one may keep all target bits separately and attach a variable to each certified target. In either convention the scalar marker weight multiplies the classical entry $t^{|\pi_\mu\vee\pi_\nu|}$ for a compatible pair. In the product basis the weighted matrix consequently has the form $G^N_{U,\boldsymbol z}\otimes F_U$.

\subsection{A four-terminal rule retaining partial minor models}
\label{sec:four-terminal-models}

A partition of four terminals need not record whether two prescribed
paths can be realized disjointly. We give a finite rule for each fixed
simple target $P$. It uses more data than
Proposition~\ref{prop:torso-marker}, but makes no connectivity
assumption on $P$. The terminals are vertices of the state graph on which the minor is tested. We describe the construction for a boundary $U$ of
size $r$. The example below takes $r=4$.

Adler et al.~\cite[pp.~7021--7023]{adler} use rooted packings to track
which branch sets meet a boundary, their connectivity, and which target
adjacencies have already been realized. Here the set of realizable
records is attached to one spanning edge state. The test on the glued
records selects its marked coefficient, with the Negami weight and the
specified ribbon evaluation.

Let $L$ be the underlying graph of a fixed edge state. Assign each
vertex of $L$ a label in $V(P)\cup\{*\}$, with $*$ meaning unused.
For $i\in V(P)$ write $L_i$ for the subgraph induced by all vertices
labelled $i$, and $U_i$ for the boundary vertices with this label.
Record the following information about the assignment:
\begin{enumerate}[label=(\roman*)]
\item the labels on $U$ and, for each $i$, the partition $\pi_i$ of
      $U_i$ induced by the connected components of $L_i$;
\item the value $h_i=\min(2,n_i)$, where $n_i$ counts components of
      $L_i$ disjoint from $U$; thus $h_i=2$ means at least two;
\item the set $D\subseteq E(P)$ of edges $ij$ witnessed by an edge
      of $L$ with endpoint labels $i$ and $j$.
\end{enumerate}
For $U_i=\varnothing$ take the unique empty partition, with zero
blocks. Denote the set of all records realized by assignments on $L$
by $\Sigma_P(L;U)$. The partitions $\pi_i$ describe connections within one prospective branch set each. The ordinary boundary partition $\pi_A$ still records connectivity of the whole edge state for its Negami weight. The set $\Sigma_P(L;U)$, taken once per graph state, is the
additional datum carried by its splitting vector.

\begin{proposition}[Gluing partial branch-set records]
\label{prop:partial-branch-records}
Suppose $L$ and $M$ have disjoint edge sets and intersect exactly in
their common boundary $U$. To compute $\Sigma_P(L\cup_U M;U)$, take
each pair of records from the two sets that agree on all boundary
labels. Keep those labels and form the output record
\[
 \pi_i=\pi_i^L\vee\pi_i^M,\qquad
 h_i=\min(2,h_i^L+h_i^M),\qquad D=D^L\cup D^M .
\]
Collect the distinct output records. A $P$ minor occurs in the glued
state if and only if some output satisfies
\[
 D=E(P),\qquad h_i+|\pi_i|=1\quad\text{for every }i\in V(P).
\]
Thus testing the records for each target and taking the disjunction
within each finite catalogue gives an exact marker rule in
Theorem~\ref{thm:q-gluing}, in particular at a four-terminal cut.
\end{proposition}

\begin{proof}
Assignments that agree on $U$ unite to one assignment on the union,
and every global assignment restricts to such a pair. For each label
$i$, a component avoiding $U$ cannot acquire vertices from the other
piece. Connections among the remaining components are exactly the
equivalence relation generated by $\pi_i^L$ and $\pi_i^M$. Counting
the former components up to two commutes with addition followed by
truncation at two. Every edge witnessing a member of $D$ lies in one
of the two pieces. This proves equality of the computed and realized
sets of records.

The equality $h_i+|\pi_i|=1$ says that the vertices labelled $i$ form
one nonempty connected subgraph. It rejects an absent label as well
as a disconnected branch set. Different labels occupy disjoint vertex
sets, and $D=E(P)$ supplies all required adjacencies. These are exactly
the conditions for a minor model. Conversely, label the vertices of
any minor model by their branch sets and all other vertices by $*$.
Taking the induced graph on a branch set preserves its connectivity.
Its record passes the test.
\end{proof}

\paragraph{How many records are possible?}
Let $k=|V(P)|$, $e=|E(P)|$, and $r=|U|$. First count only the
boundary labels and the partitions they induce. The boundary points are distinguishable. Number the target vertices
$1,\ldots,k$. If $n_*$ terminals are unused and $n_i$ have label $i$,
the multinomial coefficient counts assignments of these labels to the
terminals. Once the labels are fixed, each $U_i$ can be
partitioned in $B_{n_i}$ ways, where $B_n$ is the $n$th Bell number and
$B_0=1$. Unused terminals have no branch-set partition. Thus the exact
number of formal boundary descriptors is
\[
 A_{r,k}=
 \sum_{\substack{n_*+n_1+\cdots+n_k=r\\ n_*,n_1,\ldots,n_k\ge0}}
 \binom{r}{n_*,n_1,\ldots,n_k}\prod_{i=1}^k B_{n_i}.
\]
For fixed labels the tuple of partitions injects into the partitions of
$U$: take the union of all its blocks and add each unused terminal as
a singleton. Each terminal has $k+1$ possible labels, including $*$.
Summing over these $(k+1)^r$ label assignments gives $A_{r,k}\le(k+1)^r B_r$.
An equivalent count first selects the $s$ used terminals, partitions
them into $b$ nonempty blocks, and assigns one of the $k$ target labels
to each block. Blocks with the same label are precisely the blocks
of its partition. Writing $S(s,b)$ for a Stirling number of the second
kind, this gives $A_{r,k}=\sum_{s=0}^r\binom{r}{s}\sum_{b=0}^s S(s,b)k^b$.
In particular,
\[
 A_{4,k}=1+15k+25k^2+10k^3+k^4,\qquad A_{4,6}=4{,}447 .
\]
The entries $h_i=0,1,2$ mean zero, one, or at least two components
of label $i$ missing the boundary. There are $3^k$ formal vectors of
these counts. For each of the $e$ target edges, one bit says whether
its adjacency has been witnessed, giving $2^e$ possible subsets $D$.
Consequently the number of individual records is bounded by
\[
 N_{r,P}= A_{r,k}\,3^k 2^e \;\le\; (k+1)^r B_r\,3^k 2^e .
\]
The product counts all formal choices of these fields. Realized
records can satisfy additional constraints linking the fields. The signature $\Sigma_P(L;U)$ is a \emph{set}
of records, so $2^{N_{r,P}}$ bounds the number of different signatures.
Assignments with the same data contribute a single element to this set.

The output signature is an input for further gluing. To keep only
$W\subseteq U$ exposed, retain the labels on $W$ and replace each
$\pi_i$ by its restriction to $U_i\cap W$, dropping empty blocks.
Add the number of old blocks disjoint from $W$ to $h_i$ and truncate
the result at two. The old internal components stay internal and each
lost block supplies exactly one additional internal component. The set
$D$ is unchanged. For $W=\varnothing$ this reduces the connectivity
test to $h_i=1$ for every $i$. Duplicate records after gluing or
forgetting are retained only once. A state contributes its weight to
a sector once, however many of its assignments pass the minor test.

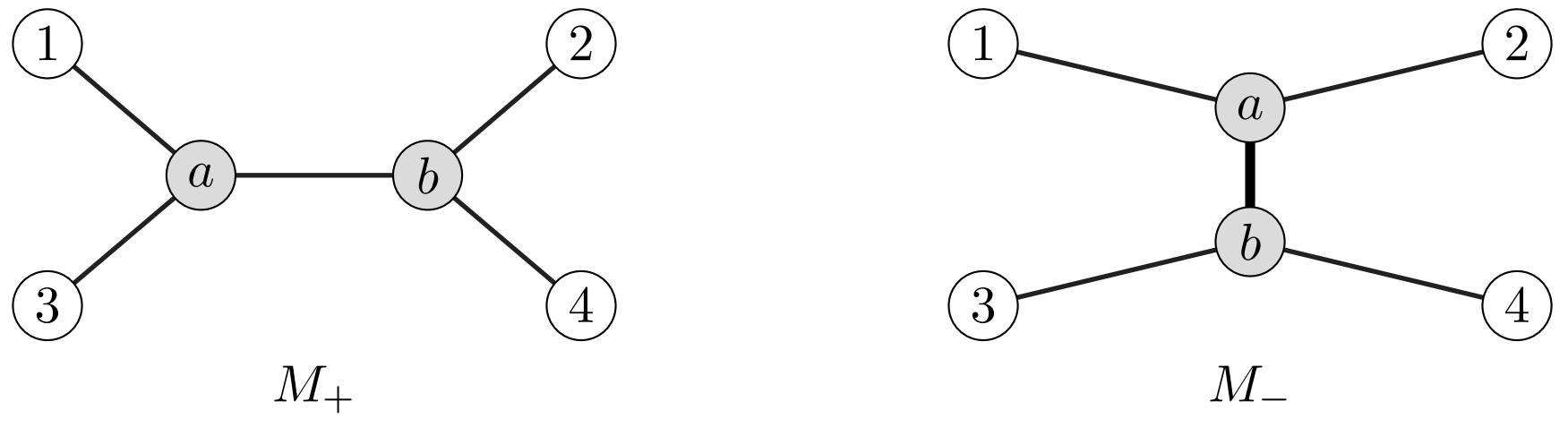
\begin{figure}[htbp]
\centering
\begin{tikzpicture}[x=1.08cm,y=.86cm,line cap=round,
    path/.style={draw=black!87,line width=1pt},
    terminal/.style={circle,draw=black,fill=white,minimum size=5.1mm,
                     inner sep=0pt,font=\normalsize},
    inner/.style={circle,draw=black,fill=black!14,minimum size=5.1mm,
                  inner sep=0pt,font=\normalsize}]
\foreach \panel/\off in {1/0,2/6.4}{
  \begin{scope}[shift={(\off,0)}]
  \coordinate (one) at (0,2.25); \coordinate (three) at (0,0);
  \coordinate (two) at (3.65,2.25); \coordinate (four) at (3.65,0);
  \ifnum\panel=1
    \coordinate (a) at (1.05,1.12); \coordinate (b) at (2.60,1.12);
    \draw[path] (one)--(a)--(three) (two)--(b)--(four) (a)--(b);
    \node at (1.825,-.72) {$M_+$};
  \else
    \coordinate (a) at (1.825,1.7); \coordinate (b) at (1.825,.55);
    \draw[path] (one)--(a)--(two) (three)--(b)--(four);
    \draw[line width=2.1pt] (a)--(b);
    \node at (1.825,-.72) {$M_-$};
  \fi
  \node[terminal] at (one) {$1$}; \node[terminal] at (two) {$2$};
  \node[terminal] at (three) {$3$}; \node[terminal] at (four) {$4$};
  \node[inner] at (a) {$a$}; \node[inner] at (b) {$b$};
  \end{scope}}
\end{tikzpicture}
\caption{Four-terminal trees glued to the same left piece
$K_6-\{13,24\}$. The white vertices are the terminals. In $M_+$ the
paths $1$--$a$--$3$ and $2$--$b$--$4$ are disjoint. In $M_-$ the paths
$1$--$a$--$b$--$3$ and $2$--$a$--$b$--$4$ share the bold edge $ab$.}
\label{fig:four-terminal-hosts}
\end{figure}

\paragraph{Two boundary states distinguished by a linkage.}
Take $U=\{1,2,3,4\}$, let $L=K_6-\{13,24\}$ on vertices
$1,\ldots,6$, and choose $P=K_6$. The right piece is a tree with two
internal vertices $a,b$ joined by an edge. In $M_+$ the leaves at $a$
are $1,3$ and those at $b$ are $2,4$. In $M_-$ the leaves at $a$ are
$1,2$ and those at $b$ are $3,4$. Consider the fully selected states.
Both trees have five edges, boundary partition $1234$, no interior
component, and no local $K_6$ minor. Even adding all six edges on the
terminals yields no $K_6$ minor on the right: the resulting six-vertex
graph is not complete. The left piece also has no local $K_6$ minor. Figure~\ref{fig:four-terminal-hosts} shows the two right pieces.

In $L\cup_U M_+$ the disjoint paths $1a3$ and $2b4$ supply the two
missing edges. For example, assign the core vertices their own labels,
give $a$ label $1$ and $b$ label $2$. This gives a passing record.
To see why $L\cup_U M_-$ fails, suppose it has six branch sets for
$K_6$. A branch set contained in $\{a,b\}$ has at most four outside
neighbors, so it cannot touch five other branch sets. Hence each branch
set contains exactly one core vertex. Relabel the target by these
vertices. Adjacency between the branch sets containing $1$ and $3$
supplies a $1$--$3$ path with all internal vertices in $\{a,b\}$;
similarly there must be a $2$--$4$ path. The two paths use disjoint
branch sets, so they are vertex-disjoint. In $M_-$ both unique paths
traverse $ab$, which is impossible. With tensor-unit decoration the
two hosts each have $13+5=18$ edges. With the singleton catalogues
$\{K_6\}$ and $\{K_7\}$, their full-edge contributions are
$z_{10}tX^{18}$ for $L\cup_U M_+$ and $z_{00}tX^{18}$ for $L\cup_U M_-$.
Neither host has a $K_7$ minor, since it has fewer than $21$ edges. The two right pieces therefore have identical ordinary partitions, local occurrence bits, and clique-completion bits before gluing.

\section{A termwise check on the sixteen states of the four-cycle}
Let the vertices of $C_4$ be $v_1,v_2,v_3,v_4$, and write its edges as
\[
 e_1=v_1v_2,\qquad e_2=v_2v_3,\qquad e_3=v_3v_4,\qquad e_4=v_4v_1.
\]
Cut along $U=\{v_1,v_3\}$, taking $e_1,e_2$ in $H$ and $e_3,e_4$ in $K$.
Encode the row subset $A\subseteq\{e_1,e_2\}$ by
$(a_1,a_2)$, where $a_i=1$ if $e_i$ is selected and $a_i=0$ otherwise.
Likewise, encode the column subset $B\subseteq\{e_3,e_4\}$ by $(b_3,b_4)$.
The labels in the tables are ordered pairs. For example, $\edgepair10$
selects $e_1$ in a row and $e_3$ in a column. In each grid the four
pairs occur in the same order on both axes.
\begin{center}
\begingroup
\renewcommand{\arraystretch}{1.2}
\setlength{\tabcolsep}{8pt}
\begin{tabular}{@{}ccccc@{}}
\toprule
\textbf{Bit pair} & \textbf{Edges in $H$} & \textbf{Edges in $K$}
    & $k_U^0$ & \textbf{Partition of $U$}\\
\midrule
$\edgepair00$ & $\varnothing$ & $\varnothing$ & 1 & $\{v_1\}\;\vert\;\{v_3\}$\\
$\edgepair10$ & $\{e_1\}$ & $\{e_3\}$ & 0 & $\{v_1\}\;\vert\;\{v_3\}$\\
$\edgepair01$ & $\{e_2\}$ & $\{e_4\}$ & 0 & $\{v_1\}\;\vert\;\{v_3\}$\\
$\edgepair11$ & $\{e_1,e_2\}$ & $\{e_3,e_4\}$ & 0 & $\{v_1,v_3\}$\\
\bottomrule
\end{tabular}
\endgroup
\end{center}
In the component table the boxed entry uses row $\edgepair10$ and
column $\edgepair10$. In set notation $A=\{e_1\}$ and $B=\{e_3\}$. Its
value is $2$, the number of components of $\{v_1v_2,v_3v_4\}$.
\begin{center}
\begingroup
\renewcommand{\arraystretch}{1.20}
\setlength{\tabcolsep}{10pt}
\setlength{\fboxsep}{1.5pt}
\begin{tabular}{@{}c@{\hspace{12pt}}*{4}{>{\centering\arraybackslash}m{.12\linewidth}}@{}}
\multicolumn{5}{c}{\textbf{Components of the glued four-cycle state}}\\[2pt]
\toprule
 & \multicolumn{4}{c}{Column $(b_3,b_4)$ on $(e_3,e_4)$}\\
\cmidrule(lr){2-5}
Row $(a_1,a_2)$ & $\edgepair00$ & $\edgepair10$ & $\edgepair01$ & $\edgepair11$\\
\midrule
$\edgepair00$ & $4$ & $3$ & $3$ & $2$\\
$\edgepair10$ & $3$ & \fbox{\(\mathbf{2}\)} & $2$ & $1$\\
$\edgepair01$ & $3$ & $2$ & $2$ & $1$\\
$\edgepair11$ & $2$ & $1$ & $1$ & $1$\\
\bottomrule
\end{tabular}
\endgroup
\end{center}
\begin{center}
\begingroup
\renewcommand{\arraystretch}{1.20}
\setlength{\tabcolsep}{10pt}
\setlength{\fboxsep}{1.5pt}
\begin{tabular}{@{}c@{\hspace{12pt}}*{4}{>{\centering\arraybackslash}m{.12\linewidth}}@{}}
\multicolumn{5}{c}{\textbf{Classical state weights for the four-cycle}}\\[2pt]
\toprule
 & \multicolumn{4}{c}{Column $(b_3,b_4)$ on $(e_3,e_4)$}\\
\cmidrule(lr){2-5}
Row $(a_1,a_2)$ & $\edgepair00$ & $\edgepair10$ & $\edgepair01$ & $\edgepair11$\\
\midrule
$\edgepair00$ & $t^4Y^4$ & $t^3XY^3$ & $t^3XY^3$ & $t^2X^2Y^2$\\
$\edgepair10$ & $t^3XY^3$ & $t^2X^2Y^2$ & $t^2X^2Y^2$ & $tX^3Y$\\
$\edgepair01$ & $t^3XY^3$ & $t^2X^2Y^2$ & $t^2X^2Y^2$ & $tX^3Y$\\
$\edgepair11$ & $t^2X^2Y^2$ & $tX^3Y$ & $tX^3Y$ & $tX^4$\\
\bottomrule
\end{tabular}
\endgroup
\end{center}

\begin{proposition}[Termwise agreement for $C_4$]\label{stmt:8.1}
For each of the sixteen states, the direct state-sum term agrees with
the corresponding term of the splitting pairing. Thus
\[
 Z^N(C_4)=t^4Y^4+4t^3XY^3+6t^2X^2Y^2+4tX^3Y+tX^4.
\]
With quantum decorations, the two calculations still agree entry by entry, since
each multiplies the entry for $C$ by the same factor $F_q((C_4)_C)$.
\end{proposition}

\begin{proof}
Apply Lemma~\ref{stmt:3.5} to each entry. For example, for row $\edgepair00$ and column $\edgepair11$, the vertex
$v_2$ is an isolated interior component on the left, there are no interior
components on the right, and the join $\{v_1,v_3\}$ has one block. Hence
$k=1+0+1=2$. The edge weight is $X^2Y^2$, giving the entry $t^2X^2Y^2$.
The same computation gives the other fifteen entries. Agreement of the
quantum factors follows from Lemma~\ref{lem:q-compose}. Summing all entries
gives the displayed formula.
\end{proof}

For row $\edgepair10$ and column $\edgepair10$, the selected edges are
$e_1=v_1v_2$ and $e_3=v_3v_4$. They form two components on the four
vertices. Both local partitions are $\{v_1\}\mid\{v_3\}$, and both
interior component counts are zero, so the join formula gives
$0+0+2=2$. With two selected and two deleted edges, this yields the
entry $t^2X^2Y^2$.

Every spanning subgraph of $C_4$ is planar, as are all of its minors.
Consequently no state has a Petersen-family minor or a $K_7$ minor.
In particular, for the $K_6$ and $K_7$ targets the marker is
$m(C)=(0,0)$ in all sixteen states. The four identities of Theorem~\ref{thm:q-gluing} give $Z^{\qMSNtag}(C_4;\boldsymbol z)=z_{00}Z^{\qMSNtag}(C_4)$, with zero coefficients in the three other sectors. Section~\ref{sec:square-marker} uses the very same sixteen edge states in a gluing that can change the closed marker.

\subsection{A marker created by gluing}
\label{sec:marked-k6}
Let $L=K_6-uv$ and let $M$ consist of the edge $uv$ on the two
terminals. Choose $K_6$ as the IL target and $K_7$ as the IK target.
For the fully retained state of $L$ the signature $(a_{K_6},b_{K_6},c)$
is $(0,1,1)$. For the retained edge on $M$ it is $(0,0,1)$.
Proposition~\ref{prop:torso-marker} gives $a_{K_6}(L\cup_U M)=1$.
If one further edge between two nonterminals is removed from $L$, the
two terminals remain connected and both local minor bits are still zero.
This time $b_{K_6}(L)=0$, and the same gluing gives $a_{K_6}=0$.
Indeed, a graph on six vertices has a $K_6$ minor precisely when it
already contains all fifteen edges. Thus the two left states have the
same boundary partition and the same local markers, but the completed-edge
test separates them.

With $F_q=1$ (the tensor-unit decoration), the sector polynomial is
\[
 Z^{\qMSNtag}(K_6;\boldsymbol z)
   =z_{00}Z^N(K_6)+(z_{10}-z_{00})\,tX^{15}.
\]
The correction moves exactly the all-retained term from sector $00$ to
sector $10$. With any compatible quantum decoration its contribution is
$z_{10}tX^{15}F_q((K_6)_{E(K_6)})$. The same scalar polynomial holds
for the full Petersen-family IL catalogue, since $K_6$ is its only
member on at most six vertices~\cite[Fig.~1.2]{hn14}.

The vertex gluing above is on the underlying graph. To realize the same
test across a regular ribbon cut, subdivide the edge $uv$ at two points
$s,r$ and place the cut through those points. One piece contains $K_6-uv$
and the two terminal stubs $us,rv$, and the other contains the single
edge $sr$. All four exposed half-edges are now univalent at the cut. In
the subdivided graph the marker is $1$ precisely when all seventeen
edges are retained. A failed stub interrupts the only replacement for
$uv$; if another core edge is missing, suppressing the degree-two
subdivision vertices yields a proper subgraph of $K_6$ or a subgraph
with pendant paths. Such paths cannot supply a $3$-connected minor. The
full term is $tX^{17}$ when $F_q=1$. This distinction records exactly
which edge set is being summed over.

\subsection{The same square as a boundary piece}
\label{sec:square-marker}
Glue the square from the beginning of this section to the core
$L=K_6-uv$ at its opposite vertices $u=v_1$ and $v=v_3$. Denote the
result by $G_\square$. Take the targets $K_6$ for IL and $K_7$ for IK,
and use the tensor-unit quantum decoration. Start with all fourteen
edges of $L$ retained. Its $K_6$ signature is $(0,1,1)$, while the
four-vertex square cannot contain either target, even after adding the
virtual edge $uv$. Proposition~\ref{prop:torso-marker} then reduces the
global IL marker to the terminal-connectivity bit of the square. Write $(a_1,a_2)$ for the upper-path bits on $(e_1,e_2)$
and $(b_3,b_4)$ for the lower-path bits on $(e_3,e_4)$, in the
indicated orders. A bit of $1$ selects its edge. For example, the
upper pair $(1,0)$ selects $e_1$ only, and the lower pair $(1,0)$
selects $e_3$ only. With the core retained, a state belongs to sector
$(1,0)$ exactly when the row or the column equals $(1,1)$. There are
$4+4-1=7$ such states among the sixteen. The other nine belong to
sector $(0,0)$. The global IK bit is zero throughout.

\begin{figure}[htbp]
\centering
\begin{tikzpicture}[x=1cm,y=.82cm,font=\normalsize,line width=1pt,
     lab/.style={font=\small,inner sep=2pt}]
\coordinate (u) at (-2.1,0); \coordinate (v) at (2.1,0);
\coordinate (x) at (0,.95); \coordinate (y) at (0,-.95);
\draw (u)--node[lab,above left] {$e_1$} (x)--node[lab,above right] {$e_2$} (v);
\draw (v)--node[lab,below right] {$e_3$} (y)--node[lab,below left] {$e_4$} (u);
\foreach \pt in {u,v,x,y}{\fill (\pt) circle (1.8pt);}
\node[left] at (u) {$u$}; \node[right] at (v) {$v$};
\node[above] at (x) {$v_2$}; \node[below] at (y) {$v_4$};
\end{tikzpicture}

\medskip
\begingroup
\renewcommand{\arraystretch}{1.15}
\setlength{\tabcolsep}{9pt}
\begin{tabular}{@{}lcccc@{}}
\toprule
& \multicolumn{4}{c}{Lower path $(b_3,b_4)$} \\
\cmidrule(lr){2-5}
Upper path $(a_1,a_2)$ & $(0,0)$ & $(1,0)$ & $(0,1)$ & $(1,1)$ \\
\midrule
$(0,0)$ & 0 & 0 & 0 & \textbf{1} \\
$(1,0)$ & 0 & 0 & 0 & \textbf{1} \\
$(0,1)$ & 0 & 0 & 0 & \textbf{1} \\
$(1,1)$ & \textbf{1} & \textbf{1} & \textbf{1} & \textbf{1} \\
\bottomrule
\end{tabular}
\endgroup
\caption{The $K_6$-minor bit for the sixteen square states attached
to the fully retained core $K_6-uv$. The upper route from $u$ passes
through $v_2$ to $v$ along $e_1,e_2$. The lower route from $v$ passes
through $v_4$ to $u$ along $e_3,e_4$. Each bold $1$ marks a state with at least one
intact terminal path.}
\label{fig:square-marker}
\end{figure}

For a two-edge path let $S=tY^2+2XY$ be the sum of the weights of its
three broken states and let $J=X^2$ be the weight of its intact state.
In a broken path the interior vertex is isolated only when both edges
are absent. The connected core contributes $tX^{14}$. The sixteen
entries of Figure~\ref{fig:square-marker} therefore contribute
\[
 tX^{14}\bigl[z_{00}S^2+z_{10}(2SJ+J^2)\bigr].
\]
The first summand is the contribution of the nine states in which both
paths are broken. The second contains the seven states with at least one
intact path. For example, the entry in row $\edgepair11$ and column $\edgepair00$ has
weight $z_{10}t^2X^{16}Y^2$, including the isolated vertex $v_4$.
If a core edge is missing, even the completed core has no $K_6$ minor,
so no further state contributes to the $10$ sector. Thus the coefficient
in the \emph{full} state sum is
\[
 [z_{10}]\, Z^{\qMSNtag}(G_\square;\boldsymbol z)
   =tX^{14}\bigl(X^4+4X^3Y+2tX^2Y^2\bigr).
\]

\section{Matrix form, Kronecker products, and the gluing defect}
\label{sec:matrix-form}

Choose a classical basis $\{e_\mu\}$, a quantum basis $\{\phi_a\}$,
and a basis $\{\psi_b\}$ of the opposite boundary space. Write
\[
 (G_U^N)_{\mu\nu}
   =\operatorname{comp}(\mu,\nu)\,t^{|\pi_\mu\vee\pi_\nu|},
 \qquad
 (F_U)_{ab}=\operatorname{ev}_U(\phi_a,\psi_b).
\]

\begin{proposition}[Kronecker factorization of the quantum gluing matrix]\label{stmt:9.1}
\label{prop:kronecker}
In the product basis $e_\mu\otimes\phi_a$ and the corresponding basis
on the right, the full gluing matrix is
\[
 G_U^{\qMSNtag}=G_U^N\otimes F_U.
\]
For scalar sector variables and an exact rule $g_U$, its marked version is
\[
 G^{\qMSNtag}_{U,\boldsymbol z}=G^N_{U,\boldsymbol z}\otimes F_U,
 \qquad
 (G^N_{U,\boldsymbol z})_{\mu\nu}=
 \begin{cases}
  t^{|\pi_\mu\vee\pi_\nu|}z_{g_U(\mu,\nu)},&
      \operatorname{comp}(\mu,\nu)=1,\\
  0,&\operatorname{comp}(\mu,\nu)=0.
 \end{cases}
\]
\end{proposition}

\begin{proof}
By Definition~\ref{stmt:6.5}, the entry in row $(\mu,a)$ and column $(\nu,b)$ is
\[
 (G_U^{\qMSNtag})_{(\mu,a),(\nu,b)}
    =(G_U^N)_{\mu\nu}(F_U)_{ab}.
\]
These are the entries of the Kronecker product. On each compatible
classical pair the marked pairing multiplies the whole quantum block
by the scalar $z_{g_U(\mu,\nu)}$. This changes only the classical matrix
and proves the second identity.
\end{proof}

\begin{definition}[Gluing defect]\label{stmt:9.2}
Given a matrix $\widehat G_U$ obtained by computation or proposed in
another theory, its \emph{gluing defect} is
\[
 D_U\coloneqq\widehat G_U-G_U^N\otimes F_U.
\]
Both matrices use the same indexed boundary spaces and basis conventions.
For a marked candidate use the marked reference
$G^N_{U,\boldsymbol z}\otimes F_U$ in the same definition.
\end{definition}

\begin{corollary}[Vanishing of the defect in the external tensor product model]\label{stmt:9.3}
For the $\qMSNtag$ theory with the pairing of Definition~\ref{stmt:6.5}, one has
$D_U=0$.
\end{corollary}

\begin{proof}
This follows from Proposition~\ref{prop:kronecker}.
\end{proof}

\begin{remark}\label{stmt:9.4}
The defect compares a candidate with a specified reference. A different
classical matrix $C_U$ can give $D_U\ne0$ even for a product
$C_U\otimes F_U$. To rule out every single product in the given
classical and quantum factors, compare the quantum blocks for fixed
classical pairs. Over a field, all nonzero blocks of a product are
proportional and hence have the same rank. Section~\ref{sec:cable-probe}
exhibits two different nonzero ranks.
\end{remark}

\subsection{A minor-controlled observable on a distinguished cable}\label{sec:cable-probe}

Use a prescribed two-strand Temperley--Lieb test cable in a separate
ball, disjoint from the graph diagram. Its two input and two output ports belong
to the probe, and its strands contribute no edges to the graph state sum. Work
over a generic coefficient field with loop value $d$, where
$d(d^2-1)\ne0$. Write $U_1$ for the cup--cap generator of
$\mathrm{TL}_2$, so that $U_1^2=dU_1$, and put
\[
 P=U_1/d,\qquad Q=1-P.
\]
These orthogonal idempotents project onto the trivial channel and its
complement. With the unnormalized closure trace,
$\operatorname{tr}(1)=d^2$ and $\operatorname{tr}(U_1)=d$, giving
$\operatorname{tr}(P)=1$ and $\operatorname{tr}(Q)=d^2-1$.

For each closed graph state $C$, take
$m=\mathbf 1[K_6\preceq G_C]$, computed across two terminals by
Proposition~\ref{prop:torso-marker}. Insert $T_0=1$ or $T_1=P$
according to this bit. The choice made from the closed graph state
acts on the cable's quantum channels. If its pieces carry
$x,y\in\mathrm{TL}_2$, their closed evaluation is
$\operatorname{tr}(xyT_m)$. In the common basis $(P,Q)$ the two
possible pairing blocks are
\[
 F^{(0)}=\operatorname{diag}(1,d^2-1),\qquad
 F^{(1)}=\operatorname{diag}(1,0).
\]
To realize both blocks, retain all edges of $K_6-uv$ on one side and omit or
retain the single edge $uv$ on the other. The respective global bits
are $0$ and $1$. The graph may carry the tensor-unit decoration, and
the independent test cable permits either $P$ or $Q$ for each fixed
classical state. Both classical join weights are nonzero for formal $t$.
In any single product $G_{\mathrm{cl}}\otimes F$, every nonzero
quantum block is a scalar multiple of the same $F$, even if the
classical data are refined. The two displayed blocks have ranks two
and one. The distinct ranks rule out such a product in the given
classical--quantum decomposition.

On the probe, $P$ is central in $\mathrm{TL}_2$; cyclicity of the
trace therefore allows the cut and the insertion to move past any
two-strand endomorphism. This argument applies to cuts retaining the
prescribed ports and the same graph state. Passage through other
vertex coupons requires the corresponding local relations.

The distinguished cable and the assignment $m\mapsto T_m$ are extra
inputs. Keeping the four-port Temperley--Lieb space and moving the
insertion onto the graph subjects it to additional braid slides
(Remark~\ref{stmt:18.1}). Section~\ref{sec:fox-constituent} instead obtains
its strands from spanning two-component cycle covers of the graph, using
a fixed Fox coloring rule.

\subsection{Two constituent triangles in a spatial \texorpdfstring{$K_6$}{K6}}
\label{sec:spatial-constituent}

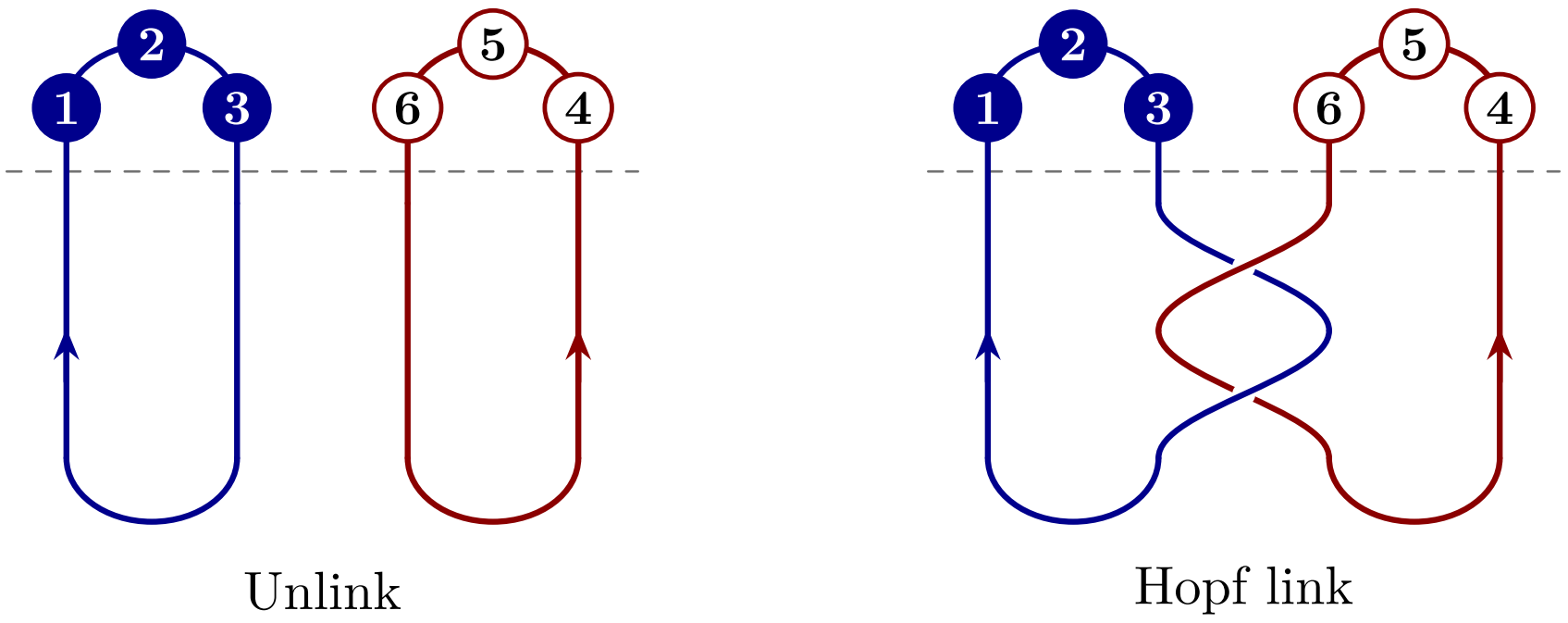
\begin{figure}[htbp]
\centering
\begin{tikzpicture}[x=1.22cm,y=.91cm,line cap=round,line join=round,
first/.style={draw=blue!55!black,line width=1.2pt},
  second/.style={draw=red!54!black,line width=1.2pt},
  firstlabel/.style={circle,draw=blue!55!black,fill=blue!55!black,
                     text=white,minimum size=4.8mm,inner sep=0pt,font=\small\bfseries},
  secondlabel/.style={circle,draw=red!54!black,fill=white,
                      text=black,line width=.9pt,minimum size=4.8mm,inner sep=0pt,
                      font=\small\bfseries}]
\foreach \panel/\off in {1/0,2/5.4}{
  \begin{scope}[shift={(\off,0)}]
  \draw[black!55,dashed,line width=.5pt] (-.35,2.25)--(3.35,2.25);
  \draw[first] (0,0)--(0,2.75) (1,2)--(1,2.75);
  \draw[second] (3,0)--(3,2.75) (2,2)--(2,2.75);
  \draw[first] (0,0) arc[start angle=180,end angle=360,radius=.5]
               (1,2.75) arc[start angle=0,end angle=180,radius=.5];
  \draw[second] (2,0) arc[start angle=180,end angle=360,radius=.5]
                (3,2.75) arc[start angle=0,end angle=180,radius=.5];
  \ifnum\panel=1
    \draw[first] (1,0)--(1,2);
    \draw[second] (2,0)--(2,2);
    \node at (1.5,-1.04) {Unlink};
  \else

  \draw[second] (2,0) .. controls (2,.38) and (1,.62) .. (1,1);
  \fill[white] (1.5,.5) circle[radius=2.35pt];
  \draw[first] (1,0) .. controls (1,.38) and (2,.62) .. (2,1);
  \draw[first] (2,1) .. controls (2,1.38) and (1,1.62) .. (1,2);
  \fill[white] (1.5,1.5) circle[radius=2.35pt];
  \draw[second] (1,1) .. controls (1,1.38) and (2,1.62) .. (2,2);
    \node at (1.5,-1.04) {Hopf link};
  \fi
  \draw[first,-{Stealth[length=2.25mm,width=1.8mm]}] (0,.6)--(0,1.02);
  \draw[second,-{Stealth[length=2.25mm,width=1.8mm]}] (3,.6)--(3,1.02);
  \node[firstlabel] at (0,2.75) {1};
  \node[firstlabel] at (.5,3.25) {2};
  \node[firstlabel] at (1,2.75) {3};
  \node[secondlabel] at (3,2.75) {4};
  \node[secondlabel] at (2.5,3.25) {5};
  \node[secondlabel] at (2,2.75) {6};
  \end{scope}}
\end{tikzpicture}
\caption{The same constituent triangles of $K_6$, oriented
$1\to2\to3\to1$ (filled vertices) and $4\to5\to6\to4$ (open
vertices), in two embeddings. In the blackboard framing, the identity
plat on the left gives the unlink, and the plat of $B_2^2$ on the right
gives the positive Hopf link. The dashed cut separates the six vertices
on the upper caps from the lower tangle. The nine edges of $K_6$ that
join the triangles carry the tensor unit and are omitted.}
\label{fig:embedded-k6-triangles}
\end{figure}

Consider the same labelled abstract graph $K_6$ in two embeddings in
$S^3$. Distinguish the disjoint triangles $123$ and $456$. In the first
embedding their union is the two-component unlink; in the second it is
the positive Hopf link. To construct them, mark three vertices on each
component of the indicated link, and use the consecutive arcs as the
six triangle edges. Draw the nine remaining edges as PL arcs, one at a
time. Each new arc can be chosen disjoint from the existing graph away
from its endpoints by general position in dimension three. Small
disjoint sectors at its endpoints accommodate the incident edges. Thus
both prescribed links extend to embeddings of the entire $K_6$. Figure~\ref{fig:embedded-k6-triangles} shows the two constituent links.
The unlinked pair in the first construction is consistent with intrinsic
linking: the Conway--Gordon assertion concerns the collection of all
disjoint pairs of cycles.

\paragraph{A coupon for every state.}
Use the bracket category over $\mathbb Q(A)$ with
$d=-A^2-A^{-2}$. Give each triangle edge the fundamental color $V$ and
each of the nine edges between the triangles the tensor-unit color.
Orient each triangle cyclically. After removing tensor-unit factors,
the coupon at a vertex is $\operatorname{id}_{\mathbf1}$ when neither
triangle edge is selected and $\operatorname{id}_V$ when both are
selected. When exactly one is selected, use the zero morphism of the
appropriate type $\mathbf1\to V$ or $V\to\mathbf1$. This convention
defines the evaluation of isolated vertices and of every partial path
with a colored end. A nonempty proper subset of the three edges of a
triangle has a vertex with exactly one selected colored edge. Its
diagram therefore evaluates to zero. For each triangle, the only
surviving choices are all three edges or none. Thus $4\cdot 2^9=2{,}048$
of the $2^{15}=32{,}768$ edge states survive this coupon test.
All retained unit edges
still count toward the classical components and the minor test.

Write $h_L$ for the unnormalized bracket of the two selected
triangles, using their blackboard framings in the standard four-strand
plat diagrams. The diagram has no self-crossings of either component.
In the bracket convention of Section~\ref{sec:tl4-check} (see also~\cite{kl}),
\[
 h_{\mathrm{unlink}}=\langle p,p\rangle=d^2,\qquad
 h_{\mathrm{Hopf}}=\langle B_2^2p,p\rangle=-d(d^2-2).
\]
The single-triangle evaluation is $d$ in either embedding. For
comparison with the oriented Jones convention, divide by $d$ and
multiply by $(-A^3)^{-w}$, with $w$ the writhe. Orient the Hopf link
with its two crossings positive. The normalized values of the unlink
and Hopf link are then $d$ and $-A^{-6}(A^4+A^{-4})$, respectively.
We use the unnormalized evaluations $h_L$ in the state sum.

\paragraph{The marked polynomial.}
Use the singleton targets $K_6$ and $K_7$ for the two minor bits. Write $Z_L^{\qMSNtag}$ for $Z^{\qMSNtag}(f(K_6);\boldsymbol z)$ with the chosen embedding $f$ and its distinguished link $L$.
Define
\[
 \begin{aligned}
 S={}&Y^6 Z^N(K_{3,3})
       +2d\,t X^3Y^3\bigl((X+Y)^3+(t-1)Y^3\bigr)^3,\\
 M={}&z_{00}\left[t^2Y^9+t\bigl((X+Y)^9-Y^9-X^9\bigr)\right]
             +z_{10}\,tX^9 .
 \end{aligned}
\]
For either embedding the entire marked sum is
\[
  Z_L^{\qMSNtag}=z_{00}S+h_L X^6 M.
\]
Indeed, with both triangles omitted, the remaining graph is $K_{3,3}$
and the six omitted edges supply $Y^6$. If exactly one triangle is
retained, contract its connected vertex set when counting components.
Each of the three opposite vertices has three edges to this connected
set. Omitting all three contributes $tY^3$, while a nonempty choice
contributes $(X+Y)^3-Y^3$. There are two choices of retained triangle
and its contribution is $d$. These cases give $S$.

When both triangles are retained, their six edges supply $X^6$.
Omitting all nine remaining edges leaves two components. Retaining at
least one joins them. On six vertices a $K_6$ minor requires all
fifteen edges, so only the fully retained state contributes to sector
$10$; sectors $01$ and $11$ are empty. This gives $M$. Subtraction isolates
the states containing both colored triangles:
\[
 Z_{\mathrm{Hopf}}^{\qMSNtag}-Z_{\mathrm{unlink}}^{\qMSNtag}
   =-d(d-1)(d+2)\,X^6 M .
\]
In particular, the two coefficients
\[
 [z_{10}\,tX^{15}]Z_L^{\qMSNtag}=h_L,\qquad
 [z_{00}\,tX^{14}Y]Z_L^{\qMSNtag}=9h_L
\]
both distinguish the embeddings at generic $A$. In the second formula
exactly one unit edge is omitted; omission of a colored edge gives
zero. Setting both occurring marker variables to $1$ also distinguishes
the ordinary decorated state sums, already in their $tX^{15}$
coefficients.

\paragraph{What the cut records.}
Take a four-punctured sphere for the plat closure and place the marked
graph vertices and the interiors of the nine unit edges in the
exterior ball. This is possible by placing the three vertices of each
component on its exterior arc before drawing the nine additional edges.
In a state with both triangles selected, the two interior tangles
induce the same endpoint partition $(12)(34)$: the braid $B_2^2$ is
pure. The labels in $(12)(34)$ are ribbon ports on the sphere. The vertices of the graph lie outside it. The two skein vectors are $p$ and $B_2^2p$, respectively. The
outside colored arcs supply the same covector, represented by $p$,
and unit-colored strands evaluate by the unit constraints. Thus the
same classical connectivity and minor sector can carry the two distinct
quantum evaluations above. This regular four-point cut computes the
quantum factor of a fixed state; edge choices in the polynomial remain
indexed by the original fifteen edges. Both embeddings use the same
specified edge colors and coupon rules. The two embeddings have the same $\operatorname{Abs}(f(G_C))$ for each
choice of $C$. In sector $10$ the coefficient $h_L$ is the framed
bracket of the two specified triangles. The marking rule and product
pairing are the same for both embeddings. The change from $p$ to
$B_2^2p$ in the boundary vector changes the evaluation of that pairing.
Thus these coefficients distinguish the two examples within the product
model of Proposition~\ref{prop:kronecker}.

\begin{corollary}[Conway--Gordon parity in full-state coefficients]
\label{cor:constituent-parity}
Let $f$ be a tame embedding of $K_6$, and let $\Lambda$ be its ten
unordered pairs of vertex-disjoint triangles. For each $\lambda\in\Lambda$,
apply the preceding coloring to that pair and write $H_\lambda(A)$ for
its coefficient $[z_{10}tX^{15}]Z^{\qMSNtag}$. Use blackboard
framing in the chosen constituent link diagram, orient its components,
and write $w_\lambda$ for its writhe. Choose $s^{1/2}=1$ at $s=1$ and normalize the coefficient by
\[
 j_\lambda(s)=(-A^3)^{-w_\lambda}H_\lambda(A)/d,\qquad s=A^{-4}.
\]
Then
\[
 \frac{1}{9}\sum_{\lambda\in\Lambda} \bigl(j_\lambda'(1)\bigr)^2
   = \sum_{\lambda\in\Lambda} \operatorname{lk}(f(\lambda))^2
   \equiv 1 \pmod 2 .
\]
\end{corollary}
\begin{proof}
Only the full-edge state has a $K_6$ certificate. The tensor-unit edges
and identity coupons reduce its evaluation to the bracket of
$f(\lambda)$, so $j_\lambda$ is its normalized Jones polynomial.
Proposition~\ref{stmt:13.2} gives $j_\lambda'(1)=-3\operatorname{lk}(f(\lambda))$
in the chosen variable. The Conway--Gordon parity theorem~\cite{cg}
says that the sum of the ten linking numbers is odd. Squaring an
integer preserves its parity. The squares also remove the choice
of component orientations.
\end{proof}
For the distinguished pair in Figure~\ref{fig:embedded-k6-triangles},
the unlink gives $j(s)=-(s^{-1/2}+s^{1/2})$ and $j'(1)=0$, whereas
the positive Hopf link gives $j(s)=-(s^{1/2}+s^{5/2})$ and $j'(1)=-3$.
If the distinguished pair is unlinked, the corollary forces a different
pair to have odd linking number. The ten evaluations use ten explicitly
specified colorings of the same spatial graph.

\section{Associativity and independence of the cut}
\label{sec:cut-independence}

\begin{proposition}[Associativity of iterated gluing]\label{stmt:10.1}
Let $H_1,H_2,H_3$ admit the same global edge states and retained marker data in both orders of gluing, with the hypotheses of Theorem~\ref{stmt:7.5} at each step. Then either order of gluing gives
the same closed evaluation under the identification
\[
 (H_1\cup H_2)\cup H_3\cong H_1\cup(H_2\cup H_3).
\]
\end{proposition}

\begin{proof}
On the classical side, the join of set partitions satisfies
\[
 (\pi_1\vee\pi_2)\vee\pi_3=\pi_1\vee(\pi_2\vee\pi_3).
\]
On the quantum side, composition satisfies
\[
 (f_3\circ f_2)\circ f_1=f_3\circ(f_2\circ f_1).
\]
Multiplication of the edge weights is also associative. The two
evaluations therefore agree term by term, and hence their state sums
agree.
\end{proof}

\begin{proposition}[Moving a cut past a braid generator]\label{stmt:10.2}
\label{prop:braid-cut}
Suppose that the cutting surface moves past one braid generator $B_i$.
If the left and right vectors transform according to
\[
 x'_L=B_i x_L,\qquad x'_R=B_i^{-1}x_R,
\]
and the pairing matrix $F_U$ satisfies
\[
 B_i^{\mathsf T}F_U=F_U B_i
 \qquad\text{or, equivalently,}\qquad B_i^{\mathsf T}F_U B_i^{-1}=F_U,
\]
then the closed pairing is invariant.
\end{proposition}

\begin{proof}
The transformed pairing is
\[
 (x'_L)^{\mathsf T}F_U x'_R
  =x_L^{\mathsf T}B_i^{\mathsf T}F_U B_i^{-1}x_R
  =x_L^{\mathsf T}F_U x_R.
\]
The two classical partitions undergo the same relabeling of the boundary
points, which preserves the number of blocks of their join.
\end{proof}

\begin{theorem}[Cut independence under local moves]\label{stmt:10.3}
\label{thm:cut-independence}
Assume that the moves identify the same admissible global edge states, with the same edge weights, component counts and any retained global markers. The edge set for the sum is fixed; moving a cut does not introduce an independently selectable subdivision edge. Assume the following quantum conditions.
\begin{enumerate}
\item Braid generators satisfy the adjoint relation of
      Proposition~\ref{prop:braid-cut}.
\item Cups and caps satisfy the snake identities of a rigid category.
\item Twists obey the chosen framing convention or writhe normalization.
\item Vertex coupons satisfy naturality and the prescribed local relations.
\end{enumerate}
Then any two cuts related by these local moves give the same closed
$\qMSNtag$ evaluation.
\end{theorem}

\begin{proof}
Decompose a move of the cut into a finite sequence of passages across
braid generators, cups and caps, twists, and coupons. Each step preserves
the pairing. For braids this follows from
Proposition~\ref{prop:braid-cut}. For cups and caps it follows from the
snake identity
\[
 (\operatorname{id}_V\otimes\operatorname{ev}_V)
 \circ(\operatorname{coev}_V\otimes\operatorname{id}_V)
 =\operatorname{id}_V.
\]
For twists it follows from the normalization, and for coupons from
naturality. A finite sequence of these moves therefore preserves the
evaluation.
\end{proof}

\subsection{A four-point check in the Temperley--Lieb category}
\label{sec:tl4-check}

Work in the Kauffman bracket category over $\mathbb Q(A)$, with the
empty diagram valued at $1$ and each contractible circle valued at
$d=-A^2-A^{-2}$. Put four fundamental endpoints in the cyclic order
$1,2,3,4$ on a disk. The two crossingless matchings
\[
 p=(12)(34),\qquad r=(14)(23)
\]
form an ordered basis of the four-point skein space. Glue two disks with
equally labeled points identified. Gluing a matching to itself produces
two circles, whereas gluing $p$ to $r$ produces one. Thus the unnormalized
bilinear closure has matrix
\[
 G=\begin{pmatrix}d^2&d\\ d&d^2\end{pmatrix},
 \qquad \langle x,y\rangle=x^{T}Gy .
\]
All coordinates below are columns, and the transpose fixes $A$.

Let $U_i$ attach a cup and a cap at positions $i,i+1$, passing the
other strands through. For instance, $U_1p=dp$, $U_1r=p$, $U_2p=r$,
and $U_2r=dr$. Consequently
\[
 U_1=U_3=\begin{pmatrix}d&1\\0&0\end{pmatrix},\qquad
 U_2=\begin{pmatrix}0&0\\1&d\end{pmatrix}.
\]
These matrices satisfy $U_i^2=dU_i$ and $U_iU_jU_i=U_i$ for
$|i-j|=1$. Take the positive generator to be
$B_i=A I+A^{-1}U_i$. The negative generator is
$B_i^{-1}=A^{-1}I+A U_i$, since the product equals
$I+(A^2+A^{-2}+d)U_i=I$. Explicitly,
\[
\begin{aligned}
 B_1=B_3&=\begin{pmatrix}-A^{-3}&A^{-1}\\0&A\end{pmatrix},&
 B_1^{-1}=B_3^{-1}&=\begin{pmatrix}-A^3&A\\0&A^{-1}\end{pmatrix},\\[3pt]
 B_2&=\begin{pmatrix}A&0\\A^{-1}&-A^{-3}\end{pmatrix},&
 B_2^{-1}&=\begin{pmatrix}A^{-1}&0\\A&-A^3\end{pmatrix}.
\end{aligned}
\]
Direct multiplication gives $U_i^{T}G=GU_i$, hence $B_i^{T}G=GB_i$
and $B_i^{T}G B_i^{-1}=G$. This verifies the cut move of
Proposition~\ref{stmt:10.2} for every pair of four-point states in these
coordinates. For a check involving two distinct matchings, take $x=p$
and $y=r$. Then
\[
 \langle B_2p,r\rangle=A d+A^{-1}d^2=-A^{-3}d,\qquad
 \langle B_2p,B_2^{-1}r\rangle=(-A^3)(-A^{-3}d)=d
 =\langle p,r\rangle .
\]

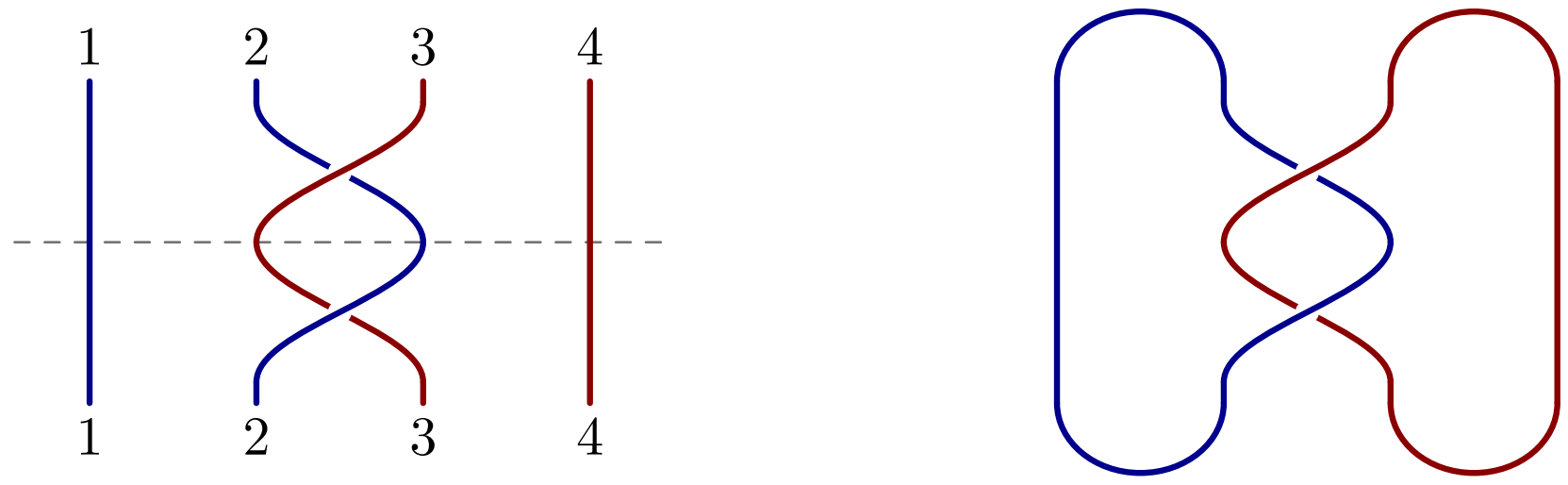
\begin{figure}[htbp]
\centering
\begin{tikzpicture}[x=1.17cm,y=.98cm,line cap=round,line join=round,
first/.style={draw=blue!55!black,line width=1.2pt},
  second/.style={draw=red!54!black,line width=1.2pt},
  endpoint/.style={inner sep=1.5pt,font=\normalsize}]
\foreach \panel/\off in {1/0,2/5.8}{
  \begin{scope}[shift={(\off,0)}]
  \ifnum\panel=1
    \draw[black!55,dashed,line width=.5pt] (-.45,1)--(3.45,1);
  \fi
  \draw[first] (0,-.15)--(0,2.15) (1,-.15)--(1,0) (1,2)--(1,2.15);
  \draw[second] (3,-.15)--(3,2.15) (2,-.15)--(2,0) (2,2)--(2,2.15);

  \draw[second] (2,0) .. controls (2,.38) and (1,.62) .. (1,1);
  \fill[white] (1.5,.5) circle[radius=2.35pt];
  \draw[first] (1,0) .. controls (1,.38) and (2,.62) .. (2,1);
  \draw[first] (2,1) .. controls (2,1.38) and (1,1.62) .. (1,2);
  \fill[white] (1.5,1.5) circle[radius=2.35pt];
  \draw[second] (1,1) .. controls (1,1.38) and (2,1.62) .. (2,2);
  \ifnum\panel=1
    \foreach \k in {1,2,3,4}{
      \node[endpoint,above] at ({\k-1},2.21) {$\k$};
      \node[endpoint,below] at ({\k-1},-.21) {$\k$};}
    \node at (1.5,-1.15) {Open braid $B_2^2$};
  \else
    \draw[first] (0,-.15) arc[start angle=180,end angle=360,radius=.5]
                 (1,2.15) arc[start angle=0,end angle=180,radius=.5];
    \draw[second] (2,-.15) arc[start angle=180,end angle=360,radius=.5]
                  (3,2.15) arc[start angle=0,end angle=180,radius=.5];
    \node at (1.5,-1.15) {Plat closure};
  \fi
  \end{scope}}
\end{tikzpicture}
\caption{On the left, positions $2$ and $3$ exchange twice. At each
crossing the strand entering from the lower left passes over the other.
The dashed line cuts the open braid into two copies of $B_2$.
On the right, the matching $p=(12)(34)$ joins $1$ to $2$ and $3$ to $4$
at both ends. The result is the blackboard-framed Hopf link.}
\label{fig:tl-plat}
\end{figure}

The braid $B_2^2$ has two successive positive crossings between
positions $2$ and $3$. The strands at positions $1$ and $4$ pass straight
through. To close the braid, join positions $1$ and $2$ by one arc and
positions $3$ and $4$ by another, using the same matching
$p=(12)(34)$ at the top and at the bottom. This four-strand plat
closure is the two-component Hopf link in Figure~\ref{fig:tl-plat}.

Expanding the two crossing operators gives
$B_2^2=A^2 I+(2+A^{-2}d)U_2$. Since
$\langle p,p\rangle=d^2$ and $\langle U_2p,p\rangle=d$, its closed
value is
\[
 \langle B_2^2p,p\rangle=A^2d^2+(2+A^{-2}d)d
 =-d(d^2-2).
\]
Cutting between the crossings gives the same expression
$\langle B_2p,B_2p\rangle$. These are framed bracket evaluations.
Dividing the closed value by $d$ gives the convention with unknot
value $1$. Passing to an oriented ambient-isotopy Jones evaluation
also requires the writhe normalization for the chosen crossing convention.
The displayed braid calculation uses a disk boundary and fundamental strands.

\paragraph{Cups, caps, twists, and a coupon.}
For the other moves in Theorem~\ref{stmt:10.3}, the same bracket convention gives explicit local checks. A cup followed by a cap closes to $d$, while the bent through-strand straightens to $\operatorname{id}_V$:
\[
 \operatorname{cap}\circ\operatorname{cup}=d,\qquad
 (\operatorname{id}_V\otimes\operatorname{cap})
 (\operatorname{cup}\otimes\operatorname{id}_V)=\operatorname{id}_V.
\]
For the curl whose $A$-smoothing creates a circle, the multiplier is $Ad+A^{-1}=-A^3$. For the opposite curl it is $A+A^{-1}d=-A^{-3}$. They are inverse multipliers. In a framed computation each is retained; in an oriented Jones computation it is cancelled by the corresponding writhe factor. Finally, for the bivalent coupon $T$ of Section~\ref{subsec:explicit-contraction}, the channel basis $(P,Q)$ gives $T=\operatorname{diag}(1,-A^4)$ as a multiplication operator. The closure form has matrix $G_{\rm ch}=\operatorname{diag}(1,d^2-1)$, so moving this coupon from one side of a pairing to the other preserves its value by $T^{\mathsf T}G_{\rm ch}=G_{\rm ch}T$. These checks establish the local identities for these particular strands and coupons.

\section{When does deletion--contraction hold?}
\label{sec:deletion-contraction}

\begin{definition}[Contraction-compatible interior edge]\label{stmt:11.1}
Let $e$ be a non-loop edge disjoint from the boundary. We call $e$
\emph{contraction-compatible} if there is a fixed linear map $C_e$ such
that, for every state $B\subseteq E(H)\setminus\{e\}$,
\[
 F_q(H_{B\cup\{e\}})=C_e F_q\bigl((H/e)_B\bigr),
\]
and the retained classical data agree under contraction. For deletion, require that omitting $e$ gives the same classical and quantum data as the corresponding state in the decorated graph $H-e$, including the prescribed coupons at any newly vacant ports. All correspondences respect admissibility. When full minor signatures are retained their agreement is a separate hypothesis.
\end{definition}

\begin{theorem}[Operator-valued deletion--contraction]\label{stmt:11.2}
\label{thm:operator-dc}
If $e$ is contraction-compatible, then
\[
 \qMSNvector{U}(H)
 =Y\qMSNvector{U}(H-e)
  +X\bigl(\operatorname{id}_{\widetilde{\mathcal N}_U}\otimes C_e\bigr)
       \qMSNvector{U}(H/e).
\]
\end{theorem}

\begin{proof}
Split the sum in Definition~\ref{stmt:6.4} according as $e\notin A$ or $e\in A$.
In the first case write $A=B$. The omitted edge contributes one extra
factor of $Y$, so these terms sum to $Y\mathbf v(H-e)$. In the second
case there is a unique expression $A=B\cup\{e\}$. Contracting a selected
non-loop edge preserves the number of connected components. The additional
selected edge contributes a factor of $X$, and contraction compatibility
identifies the quantum factor with the image under $C_e$ of the factor
for $(H/e)_B$. Adding the two contributions proves the formula.
\end{proof}

\begin{corollary}[The usual scalar recurrence]\label{stmt:11.3}
If the chosen coupon and skein conventions give
$C_e=\operatorname{id}$, then
\[
 \mathbf v(H)=Y\mathbf v(H-e)+X\mathbf v(H/e).
\]
\end{corollary}

\begin{proof}
Setting $C_e=\operatorname{id}$ in Theorem~\ref{thm:operator-dc}
makes $\operatorname{id}\otimes C_e$ the identity on the full state
space.
\end{proof}

\begin{remark}\label{stmt:11.4}
For general spatial graph vertices, contracting an edge composes its
two vertex coupons and can introduce a coefficient or a new morphism.
One must therefore determine $C_e$ and check that it is the identity
before applying the usual Negami recurrence in the quantum setting.
\end{remark}

\begin{proposition}[An edge joining two boundary points]\label{stmt:11.5}
On the ordinary partition factor, a selected edge joining boundary points $u_i$ and $u_j$ induces the block-merging operator
\[
 J_{ij}e_\pi=e_{\pi\vee\{ij\}}.
\]
If $C_e$ is the corresponding local quantum map, the operator for the
selected edge is $J_{ij}\otimes C_e$ on these two factors. Any additional marker coordinates must be updated by their own rule.
\end{proposition}

\begin{proof}
Selecting $u_i u_j$ puts $i$ and $j$ in the same connected component
regardless of the previous boundary connectivity. On partitions this
is exactly the map $\pi\mapsto\pi\vee\{ij\}$. On the quantum side,
composition with the local edge or coupon acts by $C_e$.
\end{proof}

\subsection{A non-scalar contraction operator on a three-edge path}
\label{subsec:explicit-contraction}

The four-point calculation in Section~\ref{sec:cut-independence} concerns moving
a cut past a crossing. Here the same two-strand skein algebra controls a
different operation, the contraction of an interior graph edge. The example
specifies the vertex coupon after contraction and the boundary value of
every omitted edge.

Work in the generic Kauffman bracket category over $\mathbb Q(A)$ with loop
value $d=-A^2-A^{-2}$. Each abstract graph edge carries the two-strand
cable $W=V\otimes V$, where $V$ is a single strand. An edge is selected
or omitted as a whole. Write $I=\operatorname{id}_W$ and let
$E\in\operatorname{End}(W)=\mathrm{TL}_2$ be the cup--cap diagram, so
$E^2=dE$. The two channel projectors are
\[
 P=E/d,\qquad Q=I-P.
\]
For the positive crossing take $B=A I+A^{-1}E$.
Then $BP=-A^{-3}P$ and $BQ=AQ$. Give each bivalent vertex the coupon
\[
 T=-A^3B=P-A^4Q.
\]
In particular, $TP=PT=P$. The two channels have colors $0$ and $2$,
respectively.

\paragraph{Edge and boundary conventions.}
Direct the path from its input $u$ to its output $v$, and compose maps
from right to left. An occupied edge carries the identity of $W$. For
an omitted edge the two ends terminate separately. At its source end
use $\rho=d^{-1}\operatorname{cap}:W\to\mathbf 1$, and at its target
end use $\iota=\operatorname{cup}:\mathbf 1\to W$. These maps satisfy
$\rho\iota=\operatorname{id}_{\mathbf 1}$ and $\iota\rho=P$. The
resulting disconnected diagram contributes the map $P$ from the input
port of the edge to its output port, even when one of those ports is on
the boundary of the path. In particular, all states have the same
boundary space $\operatorname{Hom}(W,W)$, naturally identified by rigidity
with a four-strand invariant space. At an omitted incoming port a coupon
is precomposed with $\iota$, and at an omitted outgoing port it is
postcomposed with $\rho$. An isolated bivalent vertex therefore evaluates
to $\rho T\iota=1$. These local rules compose in the framed ribbon
category. For example, splitting an inactive edge at an auxiliary cut
and matching its two inactive segments gives $P^2=P$. The chosen
terminations affect only the quantum diagram. Classical connectivity is
always that of the selected abstract edges.

\paragraph{Selection, deletion, and contraction.}
Let $H$ be the path with successive vertices $u,p,q,v$ and edges
$\ell=up$, $e=pq$, and $r=qv$, with $U=\{u,v\}$. Its forest states have
trivial minor markers, which we suppress. For $S\subseteq\{\ell,r\}$,
write $D_j=I$ if $j\in S$ and $D_j=P$ otherwise. In $H-e$ retain the
coupons at $p,q$ with their missing ports terminated as prescribed
above. In $H/e$ give the single merged bivalent vertex the same coupon
$T$ as every other bivalent vertex. With $F$ denoting the bracket
evaluation of the decorated state, we obtain
\begin{align*}
 F(H_{S\cup\{e\}})&=D_r T^2D_\ell,
 &
 F((H/e)_S)&=D_r T D_\ell,
 \\
 F(H_S)&=F((H-e)_S)=D_r T P T D_\ell=P .
\end{align*}
Because both possible $D_j$ commute with $T$, the fixed operator
\[
 C_e:\operatorname{End}(W)\longrightarrow\operatorname{End}(W),
 \qquad f\longmapsto T f
\]
satisfies $F(H_{S\cup\{e\}})=C_e F((H/e)_S)$ for every $S$.
The operator has eigenvalues $1$ and $-A^4$ on $P$ and $Q$.
The contracted states attain $P$ and $T$, which span
$\operatorname{End}(W)$ over $\mathbb Q(A)$. Thus both eigenvalues occur on
the space spanned by the states in this example.

\paragraph{The eight states.}
Write $e_s=e_{u\mid v}$ and $e_j=e_{uv}$ for the separated and joined
boundary partitions. The entries below are the full contributions to
$\qMSNvector{U}(H)$, including the edge weights and internal
component counts.
\begin{center}
\begin{tabular}{@{}lll@{}}
\toprule
$S$ & $e$ omitted & $e$ selected \\
\midrule
$\varnothing$ & $t^2Y^3 e_s\otimes P$ & $tXY^2 e_s\otimes P$ \\
$\{\ell\}$     & $tXY^2 e_s\otimes P$ & $X^2Y e_s\otimes P$ \\
$\{r\}$        & $tXY^2 e_s\otimes P$ & $X^2Y e_s\otimes P$ \\
$\{\ell,r\}$   & $X^2Y e_s\otimes P$ & $X^3 e_j\otimes T^2$ \\
\bottomrule
\end{tabular}
\end{center}
For example, selecting $e$ alone leaves the component on $p,q$ disjoint
from the boundary, giving the factor $t$ in the first row. The two
boundary vertices belong to one component precisely when all three
edges are selected. Abbreviating the splitting vectors to $\mathbf v$,
the table and the two smaller paths or pieces give
\begin{align*}
 \mathbf v(H-e)&=(tY+X)^2 e_s\otimes P,\\
 \mathbf v(H/e)&=(tY^2+2XY)e_s\otimes P+X^2e_j\otimes T,\\
 \mathbf v(H)&=(t^2Y^3+3tXY^2+3X^2Y)e_s\otimes P
                 +X^3 e_j\otimes T^2 .
\end{align*}
Using $C_e(P)=P$ and $C_e(T)=T^2$ gives
$\mathbf v(H)=Y\mathbf v(H-e)+X(\operatorname{id}\otimes C_e)\mathbf v(H/e)$,
as in Theorem~\ref{thm:operator-dc}. In the joined channel replacing $C_e$
by the identity would replace $T^2$ by $T$, which differs at generic $A$.
This computation concerns framed graphs with the displayed coupons
and the rule assigning $T$ to a merged vertex. A convention that retains
the composite coupon $T^2$ on that vertex absorbs the operator into the
decoration of the contracted graph.

\section{Forgetful maps and recovery of known invariants}
\label{sec:forgetful-maps}

\begin{definition}[Classical forgetful map]\label{stmt:12.1}
On the universal module with quantum diagrams as a formal basis, let
$\varepsilon_q$ be the linear map that sends each diagram to $1$.
For this map to descend to the skein quotient, we further require
that every skein relation lie in its kernel.
\end{definition}

\begin{proposition}[Recovery of the Negami state sum]\label{stmt:12.2}
Extend the classical partition module to $R$, and define $p_*\colon\widetilde{\mathcal N}_U\to R\otimes_{R_0}\mathcal N_U$ by $p_*(e_\mu)=e_{\pi_\mu}$. For the all-subsets model, if $\varepsilon_q(F_q(H_A))=1$ for every state, then
\[
 (p_*\otimes\varepsilon_q)\qMSNvector{U}(H)
 =\mathbf v_U^N(H).
\]
In particular, when $U=\varnothing$,
\[
 \varepsilon_q\bigl(Z^{\qMSNtag}(G)\bigr)
 =\sum_{A\subseteq E(G)}t^{k(A)}X^{|A|}Y^{|E|-|A|},
\]
which recovers the Negami-type state sum in our conventions.
\end{proposition}

\begin{proof}
Sending the quantum factor in each term of Definition~\ref{stmt:6.4} to $1$ and projecting the complete datum to its partition leaves
the defining sum for the classical splitting vector. When the boundary
is empty there is only one partition basis element, giving the displayed
formula.
\end{proof}

\begin{remark}\label{stmt:12.3}
The existence of a map sending every skein class to $1$ requires the
hypothesis above. A direct way to implement the forgetful map is to keep
a separate channel for it in the formal diagram module before taking
the skein quotient.
\end{remark}

\section{Extracting \texorpdfstring{$\operatorname{lk}^2$ and $a_2$}{lk squared and a2} from the Jones jet}

Let $V_L(s)$ be the normalized Jones polynomial and write $s=1+\epsilon$.
The $r$-jet is the expansion obtained by substituting $1+\epsilon$ for $s$
and taking coefficients in the truncated Taylor ring
\[
 J^r \coloneqq \Bbbk[\epsilon]/(\epsilon^{r+1}).
\]

\begin{lemma}[Jets preserve gluing]\label{stmt:13.1}
Apply the coefficient-ring homomorphism
$j_r\colon f(s)\mapsto f(1+\epsilon)\bmod\epsilon^{r+1}$ entrywise.
Then
\[
 j_r(AB)=j_r(A)\,j_r(B),\qquad
 j_r(\langle x,y\rangle)
 =\langle j_r(x),j_r(y)\rangle_{j_r(F)}.
\]
\end{lemma}

\begin{proof}
The homomorphism $j_r$ preserves addition and multiplication. Matrix
multiplication and evaluation of a bilinear pairing involve only finite
sums and products, so both commute with $j_r$.
\end{proof}

\begin{proposition}[Low-order coefficients in the chosen Jones convention]\label{stmt:13.2}
Normalize the unknot by $V_{\bigcirc}(s)=1$, and choose orientations and
signs so that $[\epsilon]V(1+\epsilon)=-3$ for the positive Hopf link.
Then, for a knot $K$ and an oriented two-component link $L$,
\begin{align*}
 a_2(K)
  &= -\tfrac13[\epsilon^2]V_K(1+\epsilon)
   = -\tfrac16 V_K''(1),\\
 \operatorname{lk}(L)
  &= -\tfrac13[\epsilon]V_L(1+\epsilon)
   = -\tfrac13 V_L'(1).
\end{align*}
\end{proposition}

\begin{proof}
For a knot, Murakami's second-derivative formula~\cite[Theorem~2]{murakami}
gives $V_K''(1)=-6a_2(K)$. Taylor expansion gives
$[\epsilon^2]V_K(1+\epsilon)=V_K''(1)/2$.
For two components his first-derivative formula~\cite[Theorem~1]{murakami}
gives three times the linking number up to the choice of Jones
variable. Evaluating the positive Hopf link in Section~\ref{sec:spatial-constituent}
fixes the sign as $-3$ in our convention.
\end{proof}

\begin{proposition}[Data needed to recover $\operatorname{lk}^2$]\label{stmt:13.3}
Use a separate marker for each constituent link $\lambda$ and record its
first-order coefficient
\[
 \ell_\lambda=-\tfrac13[\epsilon]V_{f(\lambda)}(1+\epsilon).
\]
Then the diagonal quadratic form on the $\operatorname{Sym}^2$ channel
recovers
\[
 \sum_\lambda \operatorname{lk}(f(\lambda))^2=\sum_\lambda \ell_\lambda^2.
\]
\end{proposition}

\begin{proof}
For each marker, $\ell_\lambda=\operatorname{lk}(f(\lambda))$.
Take the tensor square while retaining the markers and pair only
terms with the same marker. The resulting diagonal pairing gives the
sum of squares. Squaring the aggregate $\sum_\lambda\ell_\lambda$
would introduce cross terms $2\ell_\lambda\ell_{\lambda'}$, which
explains why the markers must be kept separate.
\end{proof}

\section{Fusion matrices and uniqueness for four and six boundary points}

\subsection{Conventions for generic \texorpdfstring{$U_q(\mathfrak{sl}_2)$}{Uq(sl2)}}

Set $A=s^{-1/4}$ and write
\[
 \delta=A^2+A^{-2},\qquad g=\delta^2-1,\qquad
 \kappa=\delta^2-2=A^4+A^{-4}.
\]
All pairings in this section are bilinear over a generic coefficient field. The transpose fixes $A$ and $\delta$. With the same basis on both sides, a change of basis acts by congruence. Different choices on the two sides give the corresponding left and right basis changes. We specify the normalization in each case~\cite{kl}.

\begin{proposition}[Fusion Gram matrix for four fundamental colors]\label{stmt:14.1}
The space of four fundamental colors with total color $0$ has a basis
$\Phi_0,\Phi_2$ of intermediate colors $0,2$. For the unnormalized
circle evaluation $d=-\delta$, its closure matrix in the common channel
basis specified below is $\operatorname{diag}(1,g)$. With the right
channel vectors scaled by $d^{-1}$, the gluing matrix is
\[
 F_4=-\delta^{-1}\begin{pmatrix}1&0\\0&g\end{pmatrix}.
\]
\end{proposition}

\begin{proof}
Use the matching basis $p,r$ and closure matrix $G$ from
Section~\ref{sec:tl4-check}. The decomposition
$1\otimes1\cong0\oplus2$ is realized by
\[
 \Phi_0=d^{-1}p,\qquad \Phi_2=r-d^{-1}p .
\]
Indeed, $U_1/d$ is the projector onto the color-$0$ channel, and
$U_1\Phi_0=d\Phi_0$, $U_1\Phi_2=0$. The three closures
$\langle p,p\rangle=\langle r,r\rangle=d^2$ and
$\langle p,r\rangle=d$ give
\[
 \langle\Phi_0,\Phi_0\rangle=1,\qquad
 \langle\Phi_0,\Phi_2\rangle=0,\qquad
 \langle\Phi_2,\Phi_2\rangle=d^2-1=g.
\]
Use $\Phi_i$ on the left and $d^{-1}\Phi_j$ on the right. In particular,
the bottom-right entry, in the color-$2$ channel, is
\[
 (F_4)_{2,2}=\langle\Phi_2,d^{-1}\Phi_2\rangle
             =d^{-1}\langle\Phi_2,\Phi_2\rangle
             =d^{-1}g=-\delta^{-1}g .
\]
Here the indices $0,2$ denote channel colors. Bilinearity over
$\mathbb Q(A)$ pulls $d^{-1}$ out unchanged. This yields the displayed $F_4$. Alternatively, after
adjoining $\sqrt d$ to the field, use $\Phi_i/\sqrt d$ on both sides
to obtain the same symmetric Gram matrix. Both computations use
bilinear evaluation and ordinary transpose.
\end{proof}

\begin{proposition}[Three channels for color-$2$ cables]\label{stmt:14.2}
Two color-$2$ objects fuse according to the rule
\[
 2\otimes2\cong0\oplus2\oplus4.
\]
The Jones--Wenzl projectors $P_0,P_2,P_4$ for these channels are
mutually orthogonal. With $[n]=(A^{2n}-A^{-2n})/(A^2-A^{-2})$ and the quantum-trace convention,
\[
 \operatorname{tr}_q(P_j)=\Delta_j=(-1)^j[j+1]\qquad(j=0,2,4).
\]
The Gram matrix for gluing is therefore
$\operatorname{diag}(\Delta_0,\Delta_2,\Delta_4)$.
\end{proposition}

\begin{proof}
The direct sum is the generic Clebsch--Gordan decomposition. The
projectors onto its simple summands satisfy $P_iP_j=0$ for $i\ne j$
and $P_j^2=P_j$. Closure evaluates the quantum trace, giving the
matrix entry $\operatorname{tr}_q(P_iP_j)=\delta_{ij}\Delta_j$.
\end{proof}

\begin{proposition}[Fusion Gram matrix for six fundamental colors]\label{stmt:14.3}
Over a splitting field containing the square roots used below, normalize the five fusion paths as in the proof. In the resulting basis
\[
 (\Phi_{AA},\Phi_{AB},\Phi_{BA},\Phi_{BB},\Phi_{CC}),
\]
the chosen convention gives
\[
 F_6=\operatorname{diag}(1,1,1,1,\kappa).
\]
\end{proposition}

\begin{proof}
A fusion path is specified by its sequence of simple intermediate
objects. When two trees have different paths, the projectors at the
first position where the paths differ are orthogonal. Hence the
off-diagonal entries vanish. At generic $A$ the closure pairing on this simple path space is nondegenerate, so each unscaled squared norm $n_i$ is nonzero. Take $h_1=\cdots=h_4=1$ and $h_5=\kappa$. Multiplying the $i$th path vector by $\sqrt{h_i/n_i}$ gives the displayed matrix. These five entries specify the basis normalization.
\end{proof}

\subsection{Uniqueness of the invariant quantum pairing}

\begin{theorem}[Uniqueness of the quantum Gram form for generic $q$]\label{stmt:14.4}
Work over a splitting field at generic $q$, so that the Temperley--Lieb category is split semisimple. Every nonzero bilinear form on the six-point invariant
space $\mathcal Q_6$ compatible with the Temperley--Lieb action is a
scalar multiple of $F_6$.
\end{theorem}

\begin{proof}
The six-point fusion-path space of total color~$0$ is the corresponding
simple module over the generic Temperley--Lieb algebra. Let $F$ and
$F'$ be two invariant forms with $F$ nondegenerate, and put
$T=F^{-1}F'$. Invariance gives
\[
 E_i^{\mathsf T}F=FE_i,\qquad E_i^{\mathsf T}F'=F'E_i
\]
for each generator $E_i$. Hence $TE_i=E_iT$ for every~$i$.
Schur's lemma gives $T=c\,\mathrm{id}$ on this simple module, so
$F'=cF$.
\end{proof}

\begin{remark}\label{stmt:14.5}
At a root of unity the module may cease to be semisimple, and
degenerate forms or additional invariant forms may occur.
The genericity hypothesis in Theorem~\ref{stmt:14.4} is essential.
\end{remark}

\subsection{Classification of weak defects and the strong forgetful condition}

For a set partition $\rho$, let $\operatorname{type}(\rho)$ denote
the integer partition formed by its block sizes. For $r=4,6$, define
\[
 (K_\lambda)_{\pi\sigma}
    =\mathbf{1}\bigl[\operatorname{type}(\pi\vee\sigma)=\lambda\bigr].
\]

The permutation action below relabels the classical terminals. We keep
the quantum ports in their specified order and impose Temperley--Lieb
invariance separately on the quantum form.
\begin{theorem}[Symmetric weak defects at four boundary points]\label{stmt:14.6}
Suppose the defect is invariant under all permutations of the boundary
points, depends on the classical variables only through the block
sizes of $\pi\vee\sigma$, and is Temperley--Lieb invariant in the
quantum variables. For generic $q$,
\[
 D_4=
 \bigl(c_{1111}K_{1111}+c_{211}K_{211}+c_{22}K_{22}
       +c_{31}K_{31}+c_4K_4\bigr)\otimes F_4.
\]
Thus there are five parameters for candidate weak defects.
\end{theorem}

\begin{proof}
Permuting the boundary points preserves the block-size type of the
join. The classical coefficient is constant on each of the five types
of integer partitions of~$4$, namely $1+1+1+1$, $2+1+1$, $2+2$,
$3+1$, and $4$. On the quantum side, the argument using Schur's lemma of
Theorem~\ref{stmt:14.4} gives a scalar multiple of $F_4$. Combining the two
factors gives the stated formula.
\end{proof}

\begin{theorem}[Symmetric weak defects at six boundary points]\label{stmt:14.7}
Under the same hypotheses,
\[
 D_6=\left(\sum_{\lambda\vdash6}c_\lambda K_\lambda\right)\otimes F_6,
\]
where $\lambda$ runs through the eleven block-size types in Table~\ref{tab:six-block-types}.
\end{theorem}

\begin{proof}
The block-size types are precisely the eleven integer partitions of~$6$. Permutations of the labelled terminals act transitively on individual
partitions of a given type. The hypothesis that the classical coefficient
depends only on the type of the join therefore leaves eleven possible
coefficients. The table counts the possible join partitions $\rho$ once
each. Several ordered pairs $(\pi,\sigma)$ may produce the same $\rho$. Theorem~\ref{stmt:14.4} makes the quantum factor
a scalar multiple of $F_6$.
\end{proof}

Table~\ref{tab:six-block-types} counts distinct set partitions of
the six labelled terminals. Each row specifies a tuple of block sizes
$\copyformula{λ}{\lambda}$ in nonincreasing order. Let
$\copyformula{m_j}{m_j}$ be the number of blocks of size $j$. The tuple
has $\copyformula{b = ∑_{j=1}^6 m_j}{b=\sum_{j=1}^6 m_j}$ blocks
and accounts for all six terminals, so
$\copyformula{∑_{j=1}^6 j m_j = 6}{\sum_{j=1}^6 j m_j=6}$. The formula in the table gives the multiplicity $N_\lambda$.
The numerator is $6!$. For each block of size $j$, divide by $j!$ to
forget its internal order. For each size $j$ divide once more by the
factorial of the number of blocks of that size. For example,
$(5,1)$ has two blocks and six choices for its singleton, while $(6)$
has one block and one partition. For $b=1,2,3,4,5,6$ the subtotals
are $1,31,90,65,15,1$. Their sum is the sixth Bell number,
$\copyformula{B_6 = 203}{B_6=203}$. In a subscript we concatenate
the block sizes, as in $\copyformula{K_{3111}}{K_{3111}}$ for $(3,1,1,1)$.

\begin{table}[htbp]
\centering
\caption{The 203 set partitions of six labelled boundary points, grouped by block sizes.}
\label{tab:six-block-types}
{\small Write $m_j$ for the number of blocks of size $j$.
For $j=1,\ldots,6$, put
$\copyformula{d_j}{d_j}=\copyformula{(j!)^(m_j)}{(j!)^{m_j}}
\cdot\copyformula{(m_j)!}{(m_j)!}$.}
\[
 \copyformula{N_}{N}_{\lambda}
 =6!\mathbin{/}(\copyformula{d_1}{d_1}\cdot\copyformula{d_2}{d_2}
 \cdot\copyformula{d_3}{d_3}\cdot\copyformula{d_4}{d_4}
 \cdot\copyformula{d_5}{d_5}\cdot\copyformula{d_6}{d_6}).
\]
\begingroup
\renewcommand{\arraystretch}{1.10}
\begin{tabular}{@{}l@{}}
\toprule
\textbf{Block sizes, block count, and set-partition count} \\
\midrule
(1, 1, 1, 1, 1, 1): 6 blocks, 1 set partition \\
(2, 1, 1, 1, 1): 5 blocks, 15 set partitions \\
(2, 2, 1, 1): 4 blocks, 45 set partitions \\
(3, 1, 1, 1): 4 blocks, 20 set partitions \\
(2, 2, 2): 3 blocks, 15 set partitions \\
(3, 2, 1): 3 blocks, 60 set partitions \\
(4, 1, 1): 3 blocks, 15 set partitions \\
(3, 3): 2 blocks, 10 set partitions \\
(4, 2): 2 blocks, 15 set partitions \\
(5, 1): 2 blocks, 6 set partitions \\
(6): 1 block, 1 set partition \\
\midrule
\textbf{All types: 203 set partitions} \\
\bottomrule
\end{tabular}
\endgroup
\end{table}

\begin{definition}[Strong boundarywise forgetful condition]\label{stmt:14.8}
Suppose the quantum vacuum $v_{\mathrm{vac}}$ satisfies
$v_{\mathrm{vac}}^{\mathsf T}F_Uv_{\mathrm{vac}}\ne0$ over the coefficient field. Over a general coefficient ring require this scalar to be a non-zero-divisor.
The strong boundarywise forgetful condition requires
\[
 D_U(e_\pi\otimes v_{\mathrm{vac}},e_\sigma\otimes v_{\mathrm{vac}})=0
\]
for every pair of classical boundary partitions $\pi,\sigma$.
\end{definition}

\begin{corollary}[The strong condition forces the defects to vanish]\label{stmt:14.9}
Under the symmetry hypotheses above and the strong boundarywise
forgetful condition,
\[
 D_4=0,\qquad D_6=0.
\]
\end{corollary}

\begin{proof}
For any join type $\lambda$, choose $\pi,\sigma$ realizing it.
Evaluating the expression for the weak defect on the vacuum gives
\[
 0=c_\lambda v_{\mathrm{vac}}^{\mathsf T}F_Uv_{\mathrm{vac}}.
\]
The vacuum pairing can be cancelled, so $c_\lambda=0$. This holds for
each type.
\end{proof}

\begin{remark}[Five directions undetected by the component count at six boundary points]\label{stmt:14.10}
Differences between shapes with the same number of blocks include
\begin{gather*}
 K_{3111}-K_{2211},\qquad K_{411}-K_{321},\qquad K_{222}-K_{321},\\
 K_{51}-K_{42},\qquad K_{33}-K_{42}.
\end{gather*}
The weight $t^{|\pi\vee\sigma|}$ does not distinguish these shapes.
The differences are therefore natural candidates in the search for a
weakly coupled theory. All disappear under the strong condition of
Corollary~\ref{stmt:14.9}.
\end{remark}

\section{Applications to graph families}

\begin{proposition}[Sufficient conditions for intrinsic linking and knotting in graph families]\label{stmt:15.1}
The following statements hold.
\begin{enumerate}
\item Cycles $C_n$ and wheels $W_n$ are planar, hence neither
 intrinsically linked nor intrinsically knotted.
\item $K_n$ is intrinsically linked for $n\ge6$ and intrinsically
 knotted for $n\ge7$.
\item $K_{p,q}$ is intrinsically linked if $p,q\ge4$.
\item $K_{p,q}$ is intrinsically knotted if $p,q\ge5$.
\end{enumerate}
\end{proposition}

\begin{proof}
For (1), a planar embedding on a sphere $S^2\subset S^3$ has only
trivial constituent knots and links. For (2), the Conway--Gordon
theorems~\cite{cg} show that $K_6$ is intrinsically linked and $K_7$
is intrinsically knotted. Larger complete graphs contain them as
subgraphs. For (3), $K_{p,q}$ contains $K_{4,4}$, which in turn
contains the Petersen-family graph $K_{4,4}-e$. Intrinsic linking
follows from the Petersen-family theorem~\cite{rst}. For (4),
$K_{5,5}$ contains $K_{3,3,1,1}$ as a minor. Indeed, contracting
two independent edges produces two adjacent vertices, each adjacent
to every vertex of the remaining $K_{3,3}$. Since $K_{3,3,1,1}$ is
intrinsically knotted~\cite{foisy}, upward monotonicity under the
minor relation proves the claim.
\end{proof}

\begin{remark}\label{stmt:15.2}
Parts (3) and (4) give sufficient conditions. Part (4) guarantees
intrinsic knotting when $p,q\ge5$ and leaves open here the full
classification over all parameter values for complete bipartite graphs.
When tabulating results, use separate columns for proved sufficient
conditions, known counterexample embeddings, and undecided cases.
\end{remark}

\subsection{A consistency check for \texorpdfstring{$K_7$}{K7} using the integral Conway--Gordon identity}

Let $\Gamma_k(K_7)$ be the set of $k$-cycles in $K_7$, and let
$\Gamma_{4,3}(K_7)$ be the set of pairs consisting of a $4$-cycle
and a vertex-disjoint $3$-cycle. The integral refinement of the
Conway--Gordon theorem~\cite{nikkuni} gives, for every spatial
embedding $f$,
\begin{equation}\label{eq:integral-cg}
\begin{aligned}
 &7\!\sum_{\gamma\in\Gamma_7(K_7)}\!a_2(f(\gamma))
 -6\!\sum_{\gamma\in\Gamma_6(K_7)}\!a_2(f(\gamma))
 -2\!\sum_{\gamma\in\Gamma_5(K_7)}\!a_2(f(\gamma))\\
 &\hspace{7em}=
 2\!\sum_{\lambda\in\Gamma_{4,3}(K_7)}\!\operatorname{lk}(f(\lambda))^2-21 .
\end{aligned}
\end{equation}

\begin{proposition}[Consistency condition for a moment-curve computation]\label{stmt:15.3}
Suppose a computation for the rectilinear embedding with vertex~$i$
at $(i,i^2,i^3)$ gives
\[
 A_7=1,\qquad A_6=A_5=0,\qquad L_{4,3}=14,
\]
where $A_k=\sum_{\gamma\in\Gamma_k(K_7)}a_2(f(\gamma))$ and
$L_{4,3}=\sum_{\lambda\in\Gamma_{4,3}(K_7)}
\operatorname{lk}(f(\lambda))^2$. These values satisfy
\eqref{eq:integral-cg}.
\end{proposition}

\begin{proof}
The left-hand side is $7\cdot1-6\cdot0-2\cdot0=7$, and the right-hand
side is $2\cdot14-21=7$.
\end{proof}

\begin{remark}\label{stmt:15.4}
This proposition checks the proposed values for the moment-curve
embedding without deriving them. It provides a consistency test for
comparing independent implementations of cycle markers, Jones jets,
the quadratic channel for $\operatorname{lk}^2$, and boundary gluing.
\end{remark}

\subsection{Labelled minor-marker sectors of \texorpdfstring{$K_7$}{K7}}
\label{sec:complete-marker-census}

Take the full Petersen family as the IL catalogue and $\{K_7\}$ as the
IK catalogue. Use the tensor-unit decoration, with every selected graph
state evaluating to~$1$. Write $z_{ij}$ for the formal weight of the
closed marker $(m_{\mathrm{IL}},m_{\mathrm{IK}})=(i,j)$, and abbreviate the resulting $Z^{\qMSNtag}(G;\boldsymbol z)$ to $Z_{\mathrm{mark}}(G)$. Each state is an edge subset of the
fixed labelled host. Isolated vertices count as components.

The six-vertex specialization was computed in Section~\ref{sec:marked-k6}. Here we compute the seven-vertex sectors and show how their IK certificate behaves under gluing.

\reserveblock{220pt}
\begin{proposition}[Marked state sum for $K_7$]
\label{prop:k7-marker-census}
Let $L_7$ be the sum of the ordinary weights of the IL-positive states
of $K_7$. For $m$ selected edges, let $N_m$ count those states and
$D_m$ count the disconnected ones. Both counts vanish for $m<15$.
\begin{center}
\normalfont
\renewcommand{\arraystretch}{1.08}
\begin{tabular}{@{}r@{\qquad}r@{\qquad}r@{}}
\toprule
Selected edges $m$ & IL-positive $N_m$ & Disconnected IL-positive $D_m$\\
\midrule
15 & 217 & 7\\
16 & 987 & 0\\
17 & 1,365 & 0\\
18 & 805 & 0\\
19 & 210 & 0\\
20 & 21 & 0\\
21 & 1 & 0\\
\midrule
Total & 3,606 & 7\\
\bottomrule
\end{tabular}
\end{center}
Let $S$ be the sum of the monomials
$\copyformula{N_m}{N_m}\,\copyformula{X^m}{X^m}\,\copyformula{Y^(21-m)}{Y^{21-m}}$
over $m=15,\ldots,21$. Multiplying $S$ by $t$ assigns every state the
weight of one component. The seven disconnected states each have
fifteen edges and two components. Correcting their contribution gives
\[
 \copyformula{L_7}{L_7}=tS+7(\copyformula{t^2}{t^2}-t)\,
              \copyformula{X^15}{X^{15}}\,\copyformula{Y^6}{Y^6}.
\]
Starting with weight $z_{00}$ on every state, give the IL-positive
states weight $z_{10}$. The only further change is the full graph,
whose marker weight is $z_{11}$. Hence
\[
 \copyformula{Z_mark}{Z_{\mathrm{mark}}}(\copyformula{K_7}{K_7})
 = \copyformula{z_00}{z_{00}}\,\copyformula{Z^N}{Z^N}(\copyformula{K_7}{K_7})
  +(\copyformula{z_10}{z_{10}}-\copyformula{z_00}{z_{00}})\,\copyformula{L_7}{L_7}
  +(\copyformula{z_11}{z_{11}}-\copyformula{z_10}{z_{10}})\,t\copyformula{X^21}{X^{21}}.
\]
At $t=X=Y=1$ the numbers of states in sectors $00$, $10$, and $11$
are $2{,}093{,}546$, $3{,}605$, and $1$, respectively.
\end{proposition}

\begin{proof}
By the excluded-minor characterization~\cite{rst}, it is enough to
test the members of the Petersen family of order at most seven.
These are $K_6$, $K_{3,3,1}$, and the graph obtained from $K_6$
by replacing a triangle by a $Y$ with a new center
\cite[Fig.~1.2 and pp.~473--474]{hn14}.
A seven-vertex minor model in this host uses seven singletons, so
each seven-vertex target must occur as a spanning subgraph.
A $K_6$ model has either six singletons or one branch set of size
two and five singletons. In the latter case the size-two set is an
edge $ab$. The five singleton vertices form a clique, and each is
adjacent to $a$ or~$b$.

These descriptions give an exhaustive calculation. Index the edges
by the pairs $ij$ with $1\le i<j\le7$ and form the following sets of
labelled witness masks.
\begin{enumerate}[label=(\roman*)]
\item Choose a six-vertex clique, leaving one vertex unused
      ($7$ masks).
\item Choose the apex of $K_{3,3,1}$ and split the remaining six
      vertices into two unordered triples ($70$ masks).
\item Choose the center of the $Y$ and its three neighbors in the
      original six-vertex clique. Remove the three edges of their
      triangle and add the three edges of the $Y$ ($140$ masks).
\item Choose $ab$, take the clique on the other five vertices, and
      for each of them choose one edge to $a$ or~$b$. Include $ab$
      itself ($21\cdot2^5$ choices before removing duplicates).
\end{enumerate}
Starting from the union of these masks, repeatedly insert every
one-edge extension, recording each resulting mask once. At closure,
count masks by edge number. The witness characterization proves that
this counts precisely the IL-positive states. The first three types
give the $7+70+140=217$ states at fifteen edges. Exactly the seven
six-clique masks at that level are disconnected. Any disconnected
IL-positive state must have a component on at least six vertices,
and that component must be $K_6$ when it has order six.
The accompanying program gives the displayed counts at all seven
levels and independently checks membership for all $82{,}160$ masks
with at least fifteen edges using complementary graphs.

A $K_7$ model in $K_7$ itself again uses only singletons. The full
edge set is therefore the unique IK-catalogue-positive state, and
its IL bit is $1$. Replacing its weight $z_{10}$ by $z_{11}$ gives
the second identity. Summing the counts gives $3{,}606$ IL-positive
states out of $2^{21}$.
\end{proof}

\subsubsection*{What the gluing edge contributes to each sector}
Write $e=uv$ and $H=K_7-e$. For each $A\subseteq E(H)$, pair its
boundary state first with the empty edge state and then with the
occupied edge state. Neither state on the edge piece contains a target.
Its terminal-connectivity bit changes from $0$ to $1$. Thus for each
target $P$, Proposition~\ref{prop:torso-marker} returns $a_P(H_A)$ in
the first pairing and $b_P(H_A)$ in the second. The two closures in Theorem~\ref{thm:q-gluing}
classify the $2^{20}$ left states at $t=X=Y=1$. Each row below describes
one transition when $e$ is inserted. In each ordered pair the IL bit
comes first and the IK bit second.
\begin{center}
\renewcommand{\arraystretch}{1.12}
\begin{tabular}{@{}c@{\qquad}c@{\qquad}r@{}}
\toprule
Before insertion & After insertion & Left states\\
$m(H_A)$ & $m(H_A+e)$ & \\
\midrule
$(0,0)$ & $(0,0)$ & 1,045,663\\
$(0,0)$ & $(1,0)$ & 2,220\\
$(0,0)$ & $(1,1)$ & 0\\
\addlinespace[3pt]
$(1,0)$ & $(0,0)$ & 0\\
$(1,0)$ & $(1,0)$ & 692\\
$(1,0)$ & $(1,1)$ & 1\\
\bottomrule
\end{tabular}
\end{center}
Exactly one of the $693$ left states with marker $(1,0)$ moves to
$(1,1)$. The remaining $692$ retain $(1,0)$. Summing all transitions
first by their initial marker and then by their final marker gives
the two closure counts.
\begin{center}
\renewcommand{\arraystretch}{1.10}
\begin{tabular}{@{}c@{\qquad}r@{\qquad}r@{\qquad}r@{}}
\toprule
Marker & $e$ absent & $e$ present & All closed states\\
\midrule
$(0,0)$ & 1,047,883 & 1,045,663 & 2,093,546\\
$(1,0)$ & 693 & 2,912 & 3,605\\
$(1,1)$ & 0 & 1 & 1\\
\midrule
All markers & 1,048,576 & 1,048,576 & 2,097,152\\
\bottomrule
\end{tabular}
\end{center}
Each left state gives two closed states, one with $e$ absent and one
with $e$ present. The last column therefore gives the three marker
coefficients of $Z_{\mathrm{mark}}(K_7)$ at $t=X=Y=1$.

For an independent derivation, let $g_r$ count left states with $r$
edges that acquire an IL certificate on insertion. Edge transitivity of
$K_7$ and monotonicity of the IL bit give
\[
  g_r=\frac{(r+1)N_{r+1}-(21-r)N_r}{21}.
\]
Indeed, the two terms after division by $21$ count IL-positive states
of size $r+1$ containing $e$ and of size $r$ omitting $e$,
respectively. Put $N_m=0$ outside the range of the table in Proposition~\ref{prop:k7-marker-census}. For $r=14,15,16,17,18$ the differences are
$155,690,870,430,75$ and they vanish elsewhere. There are $693$
IL-positive states omitting $e$, of which only $H$ gains the
$K_7$ certificate. The ancillary computation also checks these counts
by propagating uncompleted and completed witness masks across this cut.

Let $G_e(t,X,Y)$ sum the ordinary weights of those occupied-edge
closures that acquire an IL certificate when $e$ is inserted. Of the
155 such closures on fifteen edges, five consist of $K_6$ and an isolated
vertex. That vertex may be any of the five vertices outside $e$. The
other 150 are connected, as are all newly IL-positive closures with
more edges. Thus
\[
\begin{aligned}
 G_e(t,X,Y)={}&(150t+5t^2)X^{15}Y^6+690tX^{16}Y^5\\
              &+870tX^{17}Y^4+430tX^{18}Y^3+75tX^{19}Y^2 .
\end{aligned}
\]
In particular, $G_e(1,1,1)=2{,}220$, the count for $(0,0)\to(1,0)$.

Define $Z_{\mathrm{local}}(K_7)$ by using $m(H_A)$ on both
closures and retaining the ordinary edge and component weights. For the
$2{,}220$ transitions from $(0,0)$ to $(1,0)$, the marker weight changes
from $z_{00}$ to $z_{10}$. The single transition from $(1,0)$ to $(1,1)$
comes from the full graph $K_7$. Its ordinary weight is $tX^{21}$, and its
marker weight changes from $z_{10}$ to $z_{11}$. Thus
\[
 \copyformula{Z_mark}{Z_{\mathrm{mark}}}(\copyformula{K_7}{K_7})
 -\copyformula{Z_local}{Z_{\mathrm{local}}}(\copyformula{K_7}{K_7})
 = (\copyformula{z_10}{z_{10}}-\copyformula{z_00}{z_{00}})\,\copyformula{G_e}{G_e}(t,X,Y)
  +(\copyformula{z_11}{z_{11}}-\copyformula{z_10}{z_{10}})\,t\copyformula{X^21}{X^{21}}.
\]

This count concerns minors of labelled abstract graphs. The moment-curve check in Proposition~\ref{stmt:15.3} concerns evaluations in a specified spatial embedding.

\paragraph{A knotted Hamiltonian cycle in the twenty- and twenty-one-edge states.}\label{par:trefoil-k7}
Fix a labelled Hamiltonian cycle $C$ of $K_7$. Color its seven edges
by the fundamental object $V$ and its fourteen chords by the tensor
unit, retaining the two minor catalogues of the census. Use the vertex coupons of Section~\ref{sec:spatial-constituent}:
a single selected $V$-end gives zero, and two selected $V$-ends give
the identity. To realize $C$ as an unknot or a positive trefoil, mark
seven vertices on a PL representative of the chosen knot and add the
chords by disjoint arcs in general position. Write $k$ for the
unnormalized bracket of $C$ in the indicated blackboard framing.

The two embeddings give the same abstract graph at every state.
For either embedding the full state lies in sector $11$. Every one-edge deletion leaves
a connected graph containing a $K_6$ subgraph and having no $K_7$
minor. This accounts for $N_{20}=21$ and $D_{20}=0$ in the tensor-unit
census of Proposition~\ref{prop:k7-marker-census}. Deleting a cycle edge
creates two vertices with a single $V$-end. Deleting a chord preserves
the entire colored knot. Hence the decorated coefficients are
\[
 [z_{11}tX^{21}]Z^{\qMSNtag}=k,\qquad
 [z_{10}tX^{20}Y]Z^{\qMSNtag}=14k .
\]
For the crossingless unknot, $k=d$. Present the positive trefoil as
the closure of the two-strand braid $B^3$. On the two channel projectors
$P,Q$, the positive crossing $B$ has eigenvalues $-A^{-3},A$,
respectively, and the closure traces are $1,d^2-1$. Thus
\[
 k=\operatorname{tr}(B^3)=A^3(d^2-1)-A^{-9}.
\]
The diagrams have writhes $0$ and $3$. Applying $(-A^3)^{-w}/d$ to
their bracket values gives the Jones polynomials $1$ and
$s+s^3-s^4$, respectively, with $s=A^{-4}$.
For the trefoil, $V''(1)=-6$, so Proposition~\ref{stmt:13.2} gives
$a_2(C)=1$, the constituent-knot weight used in the preceding
Conway--Gordon identity.

\section{Transfer matrices}
\label{sec:transfer}

\begin{theorem}[Transfer-matrix computation]\label{stmt:16.1}
\label{thm:transfer}
Let $G_n$ be a family obtained by repeatedly gluing a cell $T$ with fixed boundary type $U$. Suppose that $\mathcal S_U$ is free of finite rank and that gluing one cell induces a linear operator $M_T\colon\mathcal S_U\to\mathcal S_U$. Then
\[
  \qMSNvector{U}(G_n)=M_T^{\,n}\qMSNvector{U}(G_0).
\]
Each coordinate satisfies the finite-order linear recurrence determined by the characteristic polynomial of $M_T$.
\end{theorem}

\begin{proof}
By the gluing formula of Theorem~\ref{thm:q-gluing}, attaching one cell acts linearly on boundary states. Attaching $n$ cells composes this map $n$ times, giving $M_T^{\,n}$. By the Cayley--Hamilton theorem, $M_T$ satisfies its characteristic polynomial. Applying this polynomial identity to $M_T^{\,n}\qMSNvector{U}(G_0)$ gives the recurrence in each coordinate.
\end{proof}

\begin{corollary}[Reduction in the size of the computation]\label{stmt:16.2}
If $d=\operatorname{rank}\mathcal S_U$ is fixed, iterating a $d\times d$ matrix computes the state vectors without enumerating all $2^{\lvert E(G_n)\rvert}$ edge subsets of $G_n$.
\end{corollary}

\begin{proof}
The matrix power in Theorem~\ref{thm:transfer} can be evaluated by repeated squaring or, when applicable, by diagonalization. Constructing the state space incurs an initial cost. Once it is in place, the same computation applies at every length $n$.
\end{proof}

\subsection{A two-state transfer for the \texorpdfstring{$K_6$}{K6} certificate}
\label{sec:k6-ear-transfer}

Let $G_s$ consist of the core $K_6-uv$ and $s$ internally disjoint
two-edge $u$--$v$ paths, with the tensor-unit decoration and target $K_6$.
A spanning state has a $K_6$ minor exactly when it retains every core
edge and at least one complete path. Indeed, completing the core by
$uv$ yields the target only if all fourteen core edges remain. The path
side completed by $uv$ has treewidth at most two (use the bags
$\{u,v,w\}$ attached to a central bag $\{u,v\}$), so it has no $K_6$ minor.
Proposition~\ref{prop:torso-marker} now gives the criterion.

With the core retained, an empty path contributes $tY^2$, a one-edge
path contributes $2XY$ across its two choices, and a complete path
contributes $X^2$. Put $a=tY^2+2XY$ and $b=X^2$. The coordinates
count states with no complete path and with at least one, respectively. Thus
\[
 \binom{w_{s,0}}{w_{s,1}}
 =\begin{pmatrix}a&0\\ b&a+b\end{pmatrix}^{s}
   \binom{tX^{14}}{0},
 \qquad w_{s,1}=tX^{14}\bigl((a+b)^s-a^s\bigr).
\]
At $t=X=Y=1$, the positive sector has $4^s-3^s$ states.
An exhaustive contraction test verifies every coefficient for $s=0,1,2$.
At $s=1000$, the transfer gives a $2000$-bit count congruent to
$692{,}330{,}323$ modulo $1{,}000{,}000{,}007$, among $2^{2014}$ edge states.

\section{Implementation and directions for further work}

\label{sec:implementation}
Fix the coefficient ring, the targets, the ribbon data and the allowed
boundary occupations before enumerating states. A boundary with four
fundamental strands has fifteen classical partitions and, in the
Temperley--Lieb example, two skein basis vectors. At six fundamental
strands of total color zero the corresponding dimensions are 203 and five.
Marking adds partial-model information to the classical coordinates.

\subsection{Record bounds and comparison with decision algorithms}
Fix one target $P$ with $k$ vertices and $e$ edges, and put $r=|U|$.
Section~\ref{sec:four-terminal-models} gives the bound
\[
    N_{r,P}=A_{r,k}\,3^k2^e
         \le (k+1)^r B_r\,3^k2^e
\]
on individual records. For fixed $r$ and $P$ this bound is independent
of the size of the host graph. For $r=4$ and $P=K_6$ it gives
$4{,}447\cdot3^6\cdot2^{15}=106{,}229{,}366{,}784$.
This number counts formal possibilities. The implementation stores
only records realized on the piece. For fixed $k\ge1$, the boundary
count $A_{r,k}$ includes all $B_r$ partitions obtained by giving every
terminal the same label. Growing the target adds a bit for every target
edge and a capped count for every target vertex.

For a fixed pair of graph states, group the two record sets by the
assignment $\ell$ of labels to their common boundary. If the corresponding
group sizes are $n_H(\ell)$ and $n_K(\ell)$, a direct join examines
at most
\[
           J_P=\sum_\ell n_H(\ell)n_K(\ell)
\]
candidate pairs, with equality if all candidates are visited. Pack the
witnessed-edge set into $\lceil e/b\rceil$ words of $b$ bits, and store
one label and block representative per terminal. A record then occupies
$O(1+r+k+\lceil e/b\rceil)$ words when the indices fit in a word.
The partition join can be computed by identifying block representatives
with a disjoint-set data structure. Together with the capped additions
and bitwise union this costs
$O(1+r\alpha(2+r)+k+\lceil e/b\rceil)$ amortized word operations per
candidate, where $\alpha$ is the inverse Ackermann function. These counts
cover forming a candidate. Duplicate removal, inclusion comparisons in
pruning, polynomial arithmetic, and ribbon evaluation incur further work.

For targets $P_1,\ldots,P_d$ these joins can be run independently on the
same pair of graph states, with at most $\sum_j J_{P_j}$ candidate
pairs in total. Aggregating weighted edge states requires the tuple
of their target signatures to preserve correlations between targets,
together with the other boundary data of Definition~\ref{stmt:6.1}.
There are at most $\prod_j 2^{N_{r,P_j}}$ formal tuples of target
signatures. This bound is only for that part of the boundary space.

With $n$ vertices and $m$ edges, a direct initialization for one edge
state tries $(k+1)^n$ vertex labelings for $P$. Direct summation also
has $2^m$ edge states. Small vertex boundaries and a short target list
make sparse records a practical starting point; repeated cells can
reuse a transfer operator once constructed. At two terminals and for
the $3$-connected targets of Proposition~\ref{prop:torso-marker}, the
common connectivity bit and two bits per target suffice instead.
The labelled $K_7$ census exploits its small host and explicit witness
masks. It does not benchmark general record enumeration at large width.

\paragraph{Comparison with minor decision algorithms.}
Let $w$ denote the width of a supplied branch decomposition of a host
with $m$ edges and let $k=|V(P)|$. For $w\ge2$, Adler et
al.~\cite[Theorem~1]{adler} decide whether this host contains $P$ in time
\[
 O\!\left(w^{2w+1}\,k^{2w}\,2^{2k^2}\,m\right).
\]
Their rooted packings at a cut are Boolean. One successful combination
suffices at the root. Our $r$ measures a single interface, while $w$ bounds
every interface of the supplied decomposition. The formal bound $N_{r,P}$
counts individual records at one interface and so has different units
from this running time.

For counts at $t=X=Y=1$ with tensor-unit evaluations, a fixed finite
catalogue also falls within the bounded-treewidth counting theorem of
Courcelle, Makowsky and Rotics~\cite[Theorem~32]{cmr}. The free edge set
specifies the state. Its target branch sets are existentially quantified
vertex sets, required to be disjoint, nonempty and connected by retained
edges, with a retained edge for each target adjacency. Counting the
free edge sets satisfying the sector formula counts each state once,
irrespective of its number of witnesses. At fixed treewidth and fixed
formula their bound is linear in host size in the stated arithmetic
model. A symbolic sector polynomial has additional coefficient costs,
and its ribbon factors require the embedding and decoration.

\subsection{Pruning and four-terminal joins}
For one graph state and target $P$, group realized records by their
boundary labels, partitions $(\pi_i)$, and capped counts $(h_i)$.
Within a group keep only inclusion-maximal witnessed-edge sets $D$.
Indeed, if $D\subseteq D'$, every compatible continuation of the first
record is compatible with the second, gives the same connectivity data,
and replaces $D\cup D_{\rm next}$ by a superset. Forgetting terminals
preserves this comparison. Thus any successful continuation through
the discarded record also succeeds through a retained one; all retained
records were realized. One may also discard a record with $h_i=2$, or
with $h_i=1$ and $U_i\ne\varnothing$, since an internal component
can never meet another component of that label through a later gluing.
These reductions preserve the minor occurrence bit, including its zero
value, under subsequent gluings and forgettings. The full signature requires the unpruned set. To reconstruct a model,
store the realizing choices along the computation. The Boolean decision
selects a sector without attaching a separate weight to each witness.

For the three fully retained pieces in Section~\ref{sec:four-terminal-models},
fix the terminal labels to $(1,2,3,4)$ in that order. Each piece has two
interior vertices and hence $7^2=49$ assignments. On the left these give
11 distinct records. All 11 pass the component-count filter, and 6 remain
after the inclusion test on $D$. On each right tree the corresponding
counts are 47, 27, and 20. The two right trees have identical counts
at every stage. Their sets of records still distinguish the required
disjoint paths. Pruned gluing, unpruned gluing, and direct minor-model
search all accept $L\cup_U M_+$ and reject $L\cup_U M_-$.

The same test over all 52 representatives of the four terminal labelings
for $K_6$ gives the following counts. Permutations of this complete
target's vertex names preserve all adjacencies. Within each boundary
labeling we count the product of the numbers of distinct left and right
records, without stopping after the first successful model.
\begin{center}
\begingroup\small
\begin{tabular}{@{}l@{}}
\toprule
Right tree (leaf pairs) \textbar{} Raw, viable and pruned candidate pairs\\
\midrule
M+ (13\textbar{}24) \textbar{} 25115 raw \textbar{} 20305 viable \textbar{} 14306 pruned\\
M- (12\textbar{}34) \textbar{} 24200 raw \textbar{} 19839 viable \textbar{} 13809 pruned\\
\bottomrule
\end{tabular}
\endgroup
\end{center}
Raw counts include every compatible pair. Viable counts apply the component
filter. Pruned counts also retain only maximal witnessed-edge sets. The integer counts are given without digit grouping in the table
and in \path{anc/pruning_observed.csv}.
These are record pairs to examine for one fully selected state of
each glued graph. The $K_7$ census of Section~\ref{sec:complete-marker-census}
uses a different computation on labelled edge subsets.

\reserveblock{126pt}
\subsection{Benchmarks on stars and their variants}\label{sec:record-benchmarks}
Take $r$ labelled terminals $0,\ldots,r-1$. In the first host, each
piece is a star on these leaves, with a distinct interior centre.
Their union is $K_{2,r}$. The \emph{fan} variant adds the path
$0,1,\ldots,r-1$ to the left star. The \emph{subdivided} variant
subdivides the arm at terminal $0$ once in each star. The fan has one interior vertex per piece, and the subdivided host has two.
The ten cases below have four to eight host vertices.
Every edge is retained. Each run measures the minor test for one target
and one graph state.

For each of the $(|V(P)|+1)^r$ assignments of named target vertices
or the unused symbol to the boundary, the program enumerates the interior
assignments on both sides and deduplicates their records. It visits
every candidate pair in that boundary group, even after finding a
model. The pruned run first applies the component and inclusion
filters above to each side. Both results are checked by a separate
enumeration of connected branch sets in the whole graph.
Each row of Table~\ref{tab:record-benchmarks} reports one benchmark case from \path{anc/benchmark_observed.csv}.
The first four columns identify the boundary width $r$, the target, the
host and the outcome of the minor test. The next three column groups
compare raw and pruned candidate-pair counts, median elapsed time in
milliseconds, and peak Python allocation in KiB. Target codes K3, C4,
K4 and K5 denote $K_3,C_4,K_4$ and $K_5$, respectively.

\begin{table}[!htb]
\centering
\caption{Minor tests on stars and their variants, before and after pruning.}
\label{tab:record-benchmarks}
\smallskip
\begingroup
\small
\renewcommand{\arraystretch}{1.17}
\setlength{\tabcolsep}{4.5pt}
\newcommand{\benchgap}{\char23\relax}
\newcommand{\benchstart}[1]{\BeginAccSupp{method=escape,ActualText={#1}}}
\begin{tabular}{@{}rllc@{\hspace{13pt}}rr@{\hspace{13pt}}rr@{\hspace{13pt}}rr@{}}
\toprule
\multicolumn{4}{c}{\textbf{Instance}} &
\multicolumn{2}{c}{\textbf{Candidate pairs}} &
\multicolumn{2}{c}{\textbf{Time (ms)}} &
\multicolumn{2}{c}{\textbf{Peak (KiB)}} \\
\cmidrule(r){1-4}\cmidrule(lr){5-6}\cmidrule(lr){7-8}\cmidrule(l){9-10}
\benchstart{r Target Host Minor Raw Pruned Raw Pruned Raw Pruned}$r$\benchgap & Target\benchgap & Host\benchgap & Minor\benchgap & Raw\benchgap & Pruned\benchgap & Raw\benchgap & Pruned\benchgap & Raw\benchgap & Pruned\EndAccSupp{} \\
\midrule
\benchstart{2 K3 Stars Yes 172 142 0.4 0.5 14.1 14.7}2\benchgap & K3\benchgap & Stars\benchgap & Yes\benchgap & 172\benchgap & 142\benchgap & 0.4\benchgap & 0.5\benchgap & 14.1\benchgap & 14.7\EndAccSupp{} \\
\benchstart{3 K3 Stars Yes 835 577 2.0 2.6 15.8 16.4}3\benchgap & K3\benchgap & Stars\benchgap & Yes\benchgap & 835\benchgap & 577\benchgap & 2.0\benchgap & 2.6\benchgap & 15.8\benchgap & 16.4\EndAccSupp{} \\
\benchstart{4 K3 Stars Yes 3760 2488 11.0 12.8 16.7 17.3}4\benchgap & K3\benchgap & Stars\benchgap & Yes\benchgap & 3760\benchgap & 2488\benchgap & 11.0\benchgap & 12.8\benchgap & 16.7\benchgap & 17.3\EndAccSupp{} \\
\benchstart{5 K3 Stars Yes 15859 10519 56.2 66.3 19.1 19.7}5\benchgap & K3\benchgap & Stars\benchgap & Yes\benchgap & 15859\benchgap & 10519\benchgap & 56.2\benchgap & 66.3\benchgap & 19.1\benchgap & 19.7\EndAccSupp{} \\
\benchstart{4 C4 Stars Yes 14089 8401 31.5 39.1 18.7 19.5}4\benchgap & C4\benchgap & Stars\benchgap & Yes\benchgap & 14089\benchgap & 8401\benchgap & 31.5\benchgap & 39.1\benchgap & 18.7\benchgap & 19.5\EndAccSupp{} \\
\benchstart{4 K4 Stars No 14833 10225 28.9 36.3 17.9 19.1}4\benchgap & K4\benchgap & Stars\benchgap & No\benchgap & 14833\benchgap & 10225\benchgap & 28.9\benchgap & 36.3\benchgap & 17.9\benchgap & 19.1\EndAccSupp{} \\
\benchstart{4 K5 Stars No 45116 32596 80.5 102.3 19.8 21.4}4\benchgap & K5\benchgap & Stars\benchgap & No\benchgap & 45116\benchgap & 32596\benchgap & 80.5\benchgap & 102.3\benchgap & 19.8\benchgap & 21.4\EndAccSupp{} \\
\benchstart{4 K4 Fan Yes 11389 7717 33.7 39.5 19.4 20.3}4\benchgap & K4\benchgap & Fan\benchgap & Yes\benchgap & 11389\benchgap & 7717\benchgap & 33.7\benchgap & 39.5\benchgap & 19.4\benchgap & 20.3\EndAccSupp{} \\
\benchstart{4 K4 Subdivided No 222269 39633 179.3 195.9 35.9 31.6}4\benchgap & K4\benchgap & Subdivided\benchgap & No\benchgap & 222269\benchgap & 39633\benchgap & 179.3\benchgap & 195.9\benchgap & 35.9\benchgap & 31.6\EndAccSupp{} \\
\benchstart{4 K3 Subdivided Yes 31480 5392 53.0 53.4 25.2 25.2}4\benchgap & K3\benchgap & Subdivided\benchgap & Yes\benchgap & 31480\benchgap & 5392\benchgap & 53.0\benchgap & 53.4\benchgap & 25.2\benchgap & 25.2\EndAccSupp{} \\
\bottomrule
\end{tabular}
\endgroup
\end{table}

For the two stars and $P=K_3$, increasing $r$ from two to five raises
the raw pair count from 172 to 15859. At width four, the targets $C_4$
and $K_4$ use the same boundary labelings but record different edge
witnesses and give different minor decisions. Subdividing one arm on
each side preserves the $K_4$ decision and increases the raw pair count
from 14833 to 222269. Pruning removes about $82\%$ of the latter
pairs. Record construction and filtering still take time, and none of
these ten runs shows a lower median total time after pruning.

The ancillary program \texttt{benchmark\_minor.py} uses a single process.
Times are median elapsed times for five runs after warm-up, including
record generation, duplicate removal, any pruning, and the complete
join. Independent validation is outside the timer. Peak values are
Python allocations measured in separate \texttt{tracemalloc} runs.
One boundary-label group is processed at a time.
The table was obtained with CPython 3.12.14 on x86-64 Linux (AMD EPYC
9V74), hash seed zero, and the options \texttt{--repeat 5 --memory}.

The interface in these benchmarks is the shared set of leaves. The
two stars form $K_{2,r}$, whose treewidth is two for $r\ge2$. The
cost here can grow with the chosen interface even on this bounded-treewidth
family. Section~\ref{sec:k6-ear-transfer} uses the two centres as terminals
for a repeated cell and gives exact sector counts along a long family.
The two calculations illustrate why the choice of cut, as well as the
number of host vertices, matters in practice.

\section{Scope and nonproduct pairings on graph strands}

\begin{center}
\small
\renewcommand{\arraystretch}{1.21}
\begin{tabularx}{\linewidth}{@{}>{\raggedright\arraybackslash}p{.32\linewidth}>{\raggedright\arraybackslash}X@{}}
\toprule
\textbf{Result} & \textbf{Hypotheses or scope} \\
\midrule
Edge-disjoint gluing formula & Reindexing of compatible edge states using the Negami join and the specified ribbon functor. Marked coefficients also require a locally computable exact rule. Proposition~\ref{prop:partial-branch-records} supplies one for finite lists of simple targets on a finite vertex boundary.\\
Kronecker factorization & Ordinary and scalar-marked pairings have one common quantum form $F_U$.\\
Nonproduct cable observable & A disjoint two-strand probe and a chosen marker-dependent insertion give blocks of ranks two and one at generic parameter.\\
Finite coloring model & For a fixed $K_6$ test and Fox coloring rule, spanning pairs of graph cycles give an ambient-isotopy invariant with nonproportional coloring blocks (Proposition~\ref{prop:fox-nonproduct}).\\
Cut independence & Fixed global states and weights, and the stated local ribbon relations.\\
Deletion--contraction & Deletion-compatible decorations and a fixed $C_e$ for the selected branch; computed on a two-strand cable in Section~\ref{subsec:explicit-contraction}.\\
Uniqueness of the quantum form & Generic $q$ and split simplicity of the corresponding Temperley--Lieb module.\\
$D_4=D_6=0$ & The strong boundarywise forgetful condition. The weak conditions leave free parameters.\\
IL marker & Complete for the full Petersen family. Smaller catalogues certify only their specified targets.\\
IK marker & Positive for the chosen IK targets. For $\{K_7\}$ the bit vanishes on the intrinsically knotted graph $K_{3,3,1,1}$ (Remark~\ref{stmt:2.8}).\\
Sign in the Jones jet & Fix the Jones variable, the orientations, and the convention for the Hopf link.\\
\bottomrule
\end{tabularx}
\end{center}

\subsection{Transport of a four-port projector}
\begin{remark}[The four-port obstruction]\label{stmt:18.1}
Section~\ref{sec:cable-probe} lets a closed minor bit select an insertion on
an auxiliary cable. Its vacuum-channel projector is central in
$\mathrm{TL}_2$. Moving that insertion between two of the graph's own
strands introduces a new slide past the neighboring strand.
In the notation of Section~\ref{sec:tl4-check}, $P=U_1/d$ fails the corresponding double slide, since
\[
  [P,B_2^2]p=(1-A^{-4})(p/d-r).
\]
Over $\mathbb Q(A)$ the right side is nonzero. The next proposition
shows what happens if one tries to remember the motion of the cut
by attaching a finite label to it.
\end{remark}

\begin{proposition}[Finite return of the cut data at generic parameter]
\label{prop:finite-return}
Work over $\mathbb Q(A)$ on the four fundamental ports of total color zero.
Let $F$ be their nonsingular closure matrix. For each classical pair $c$
write its bilinear block as $FT_c$. Suppose that for each elementary
braid $B_i$ the cut can be slid repeatedly, with invertible update $s_i$
of the data. Assume that, for every $c$ and $i$, some positive number
$n_i(c)$ of these slides returns the data to $c$. If each slide preserves
the pairing on all left and right test vectors, in the sense that
\[
 B_i^{\mathsf T}F T_{s_i c}B_i^{-1}=F T_c ,
\]
then each $T_c$ is scalar. In particular, finitely many reversibly
updated labels satisfy the return hypothesis.
\end{proposition}
\begin{proof}
Put $P_i=U_i/d$ and $\eta=-A^{-4}$. The eigenvalues of $B_i$ are
$A\eta$ on $P_i$ and $A$ on $I-P_i$, so
\[
 B_i^n=A^n\bigl(I+(\eta^n-1)P_i\bigr).
\]
Since $B_i^{\mathsf T}F=FB_i$, iterating the slide to a return
gives $[T_c,B_i^{n_i(c)}]=0$. The element $\eta$ has infinite order,
so $T_c$ commutes with every $U_i$. The $d$-eigenlines of $U_1,U_2$
are spanned by $p,r$, respectively. An operator preserving both lines
and commuting with $U_1 r=p$ has equal eigenvalues on them.
Thus all the blocks are multiples of $F$.
\end{proof}

\paragraph{Changing the cut.}
For a path of invertible cut moves with transport $R$, covariance
prescribes the new operator $R^{-1}TR$, after identifying the channel
spaces. Two paths to the same cut with the same boundary data give the
same operator exactly when the corresponding return motion
commutes with $T$. The conjugacy orbit has finite size exactly when
the stabilizer has finite index in the transport group. These are
distinct tests. For example, at generic $A$ the conjugates of $P$
under powers of $B_2^2$ are all distinct. A repetition would make $P$
commute with a positive power of $B_2^2$ and hence with $U_2$,
contradicting Remark~\ref{stmt:18.1}. Fixing a single cut removes some
return motions from the test. A closed invariant that allows the cut
to vary must still agree around every admissible loop and through cups,
caps and vertices.

\paragraph{A finite orbit at a root of unity.}
At $A=1$, $d=-2$ and $\det F=12$. The matrices $S_i=B_i=I+U_i$
are involutions satisfying the braid relations. Their action on the three
lines $\mathbb Qp,\mathbb Qr,\mathbb Q(p+r)$ permutes these lines.
The $F$-orthogonal projectors onto them, in that order, are
\[
 P_p=\begin{pmatrix}1&-1/2\\0&0\end{pmatrix},\qquad
 P_r=\begin{pmatrix}0&0\\-1/2&1\end{pmatrix},\qquad
 P_{p+r}=\tfrac12\begin{pmatrix}1&1\\1&1\end{pmatrix}.
\]
In particular $P_p=P$ and $S_2$ interchanges $P_p$ and $P_{p+r}$.
Assign $FT_c=Fc$ to the orbit label $c$, with $s_i(c)=S_i^{-1}cS_i$.
These blocks satisfy the slide equation and are nonproportional.
Their sum is $\tfrac32 F$, so summing all three labels with equal
weight restores a product. A graph construction would have to supply
a label and exhibit the distinct blocks on common realizable tests
while respecting general gluings. This calculation establishes the
finite orbit and braid covariance.

\paragraph{Other representations.}
For cabled or higher-color strands the same return argument tests the
commutant of the admissible braid powers on the chosen closure space.
A reducible representation can admit several blocks. Distinct blocks
must survive closure and occur on the actual graph states. A singular
specialization requires a further check. For example, when $d=\pm1$
the four-port form $F$ annihilates $p-dr$. Its quotient by the radical
has dimension one, so the usual closure cannot distinguish two
nonproportional forms on that quotient. An approach using indecomposable
modules and a different trace would have to prove compatibility of the
trace with sewing and with the unframed isotopy moves.

\subsection{A finite coloring model on constituent cycles}
\label{sec:fox-constituent}

There is an ambient-isotopy-invariant alternative to the generic
Temperley--Lieb pairing in Remark~\ref{stmt:18.1}. Its strands are
edges of the embedded graph. Fix the abstract target $K_6$ and the
following coloring rule once for all graphs. For an embedding $f\colon G\hookrightarrow S^3$ and $C\subseteq E(G)$
let $b(C)=\mathbf 1[K_6\preceq G_C]$ and let $\mathcal L_2(C)$
be the set of spanning $2$-regular subgraphs of $G_C$ with exactly
two components. Each element is one unordered pair of disjoint cycles.
For a link $L$, write $N_1(L)$ for its number of Fox $3$-colorings and
$N_0(L)=3$ for the number assigning a single common color to all its
arcs. Define
\[
 \mathcal F(f(G))=
 \sum_{C\subseteq E(G)}w(C)\, z_{b(C)}
       \sum_{L\in\mathcal L_2(C)} N_{b(C)}(f(L)).
\]
Here $z_0,z_1$ mark the two sectors of this one-target example. We work with a fixed finite graph, including its given degree-two vertices. Distinct
cycle covers are summed once each.

\begin{proposition}[Finite coloring gluing and two nonproportional blocks]
\label{prop:fox-nonproduct}
The polynomial $\mathcal F$ is invariant under graph relabeling and
ambient isotopy of the embedding. Its coloring counts glue by finite
boundary data. In the decomposition into graph data and coloring
data described below, two attainable marker sectors give nonproportional
forms on the same space of realizable tangle states.
\end{proposition}

\begin{proof}
The colors are elements of $\mathbb F_3$. At a crossing, an over-arc
of color $b$ sends an under-arc color $a$ to $a*b=2b-a$.
The identities $a*a=a$, $(a*b)*b=a$, and
$(a*b)*c=(a*c)*(b*c)$ give bijections of colorings across the three
Reidemeister moves, fixing colors outside the move. They also preserve
the monochromatic colorings~\cite{przytycki-fox}. This calculation uses
unoriented diagrams and requires no framing. The cycle pairs and minor
bits are determined on the abstract graph. Ambient isotopy carries each
constituent link to an isotopic link, and relabeling permutes the state
and cycle-pair summands.

For an edge-disjoint geometric gluing, choose cover edges only among
those of $C$. At a vertex boundary $U$, record the cover degrees on $U$,
the partition $\pi$ induced on $U$ by the spanning cover-edge subgraph,
and its number $k^0$ of components disjoint from $U$ (capped at three).
Require degree two at every original interior vertex. On gluing, require
the two degrees to sum to two at each shared vertex. By
Lemma~\ref{stmt:3.5}, the condition for exactly two components is
$k_L^0+k_R^0+|\pi_L\vee\pi_R|=2$, also valid with the capped counts.
At a shared vertex of degree split $1+1$ match the two strand ends.
For $2+0$ or $0+2$ join the two incident cover edges on their own side.
A regular cut through original edges inherits their occupancy and weight
once, and temporary ports are not new spanning vertices. Restriction
and union identify each accepted $(C,L)$ exactly once.
Proposition~\ref{prop:partial-branch-records} supplies the marker at
vertex boundaries.

For a fixed partial cover $T$ with $r$ matched strand ports (only
$1+1$ splits at a vertex cut), use the
two summands of
\[
 W_r=\mathbb Q^{\mathbb F_3^r}\oplus\mathbb Q^{\mathbb F_3}.
\]
In the first summand the coordinate at $\alpha$ counts all colorings of
$T$ extending the port assignment $\alpha$. In the second, the
coordinate at $a$ counts the single coloring identically $a$, including
on any component inside the piece. Keep the same global color label
through successive gluings. For $r=0$ the first summand has one
coordinate and the second still has three labels, each with value one
even for an empty piece. The dot product on the summand selected by
$b(C)$ (first for $1$, second for $0$) thus counts exactly the desired
closed colorings, including when entire components lie on one side.
At a cup or cap the two incident colors agree. A braid permutes the
finite set of port tuples, so its action on color assignments has
finite image. Counts are integer coefficients and need not be bounded.

To see both forms, choose two complementary triangles of $K_6$.
Subdivide one edge of each triangle twice and the other four triangle
edges once. This gives 14 vertices and 23 edges. Every spanning
cycle cover must use both edges at each new degree-two vertex. This
forces the two subdivided triangles and saturates every original
vertex, giving a unique cover. Retain all its edges and eight of the
nine chords between the triangles. The last chord changes $b$ from
$0$ to $1$, since suppressing the degree-two subdivisions gives $K_6-e$ or
$K_6$, respectively. To see that subdivision preserves the $K_6$ test, observe that in a
minor model a degree-two subdivision vertex used alone cannot have the
five required neighboring branch sets, and otherwise it can be
suppressed along its path. Both states are connected and their classical
weights are $tX^{22}Y$ and $tX^{23}$.

Place the four vertices on the two doubly subdivided edges at the
cut sphere. In the interior ball put just the two edges joining their
respective subdivision vertices. In the exterior put the two remaining
cycle arcs and all nine chords, routed disjointly except at their
ends as in Section~\ref{sec:spatial-constituent}. On the interior choose
the matching $p=(12)(34)$ or its image under two crossings of the
middle strands. On the exterior choose either reflected test tangle, so
the same port coordinates are compared on gluing. A pure braid
preserves the endpoint pairs. The two tangle choices are independent,
and the graph data on each side are identical for both tangle choices
at fixed chord occupancy. Switching the chord leaves all coloring
vectors unchanged.

The boundary color tuples for $p$ form
$S=\{(a,a,b,b):a,b\in\mathbb F_3\}$. The crossing permutation on the
middle coordinates is $(x,y)\mapsto(y,2y-x)$. Choose the tangle on
each side with boundary support $S$ or its twice-transformed image
$S'=\{(a,2b-a,a,b):a,b\in\mathbb F_3\}$. Each tuple in its support
extends uniquely. Both supports have size nine, and their intersection
consists of the three constant tuples. On the span of these two
realizable state vectors, the respective coloring blocks are
\[
  \copyformula{M_1 = [[9, 3], [3, 9]]}{M_1=\begin{pmatrix}9&3\\3&9\end{pmatrix}},\qquad
  \copyformula{M_0 = [[3, 3], [3, 3]]}{M_0=\begin{pmatrix}3&3\\3&3\end{pmatrix}}.
\]
Their ranks over $\mathbb Q$ are two and one. The nonzero classical
weight within each block is common to all four entries. Thus the two
blocks cannot arise from a single product of a classical matrix with
one fixed coloring form in this decomposition.
\end{proof}

The entire polynomial of this host is
\[
 \begin{aligned}
 \mathcal F(f(G))={}&3z_0X^{14}
  \left[t^2Y^9+t\bigl((X+Y)^9-Y^9-X^9\bigr)\right]\\
 &{}+z_1 tX^{23} N_1(f(L)),
 \end{aligned}
\]
where $L$ is the unique cover and $N_1(f(L))$ is the corresponding
entry of $M_1$. Precisely $2^9=512$ edge states have nonzero inner
sum. With no chord the two cycles form separate components, any of
the $510$ nonempty proper chord sets joins them in sector $0$, and
the full chord set lies in sector $1$.
In particular let $f_{00}$ use straight tangles on both sides and let
$f_{01}$ twist one side. The unique cover closes to the two-component
unlink in the first embedding and to a Hopf link in the second.
Subtraction gives the embedding-sensitive coefficient
\[
  \copyformula{ℱ(f_{00}(G)) − ℱ(f_{01}(G)) = 6 z_1 t X^{23}}{\mathcal F(f_{00}(G))-\mathcal F(f_{01}(G))=6z_1tX^{23}}.
\]

The finite coloring sum uses a reducible boundary space and an explicit
global rule on all spanning two-cycle covers. It supplies a nonproduct
example for that rule on undecorated graphs. The existence of an intrinsic coupling on the generic Temperley--Lieb space under weaker slide or sewing hypotheses remains open.

\section{Conclusion}

The marked gluing formula of Theorem~\ref{thm:q-gluing} is
\[
 Z^{\qMSNtag}(H\cup_U K;\boldsymbol z)
   =\bigl\langle\qMSNvector{U}(H),
                \qMSNvector{U}(K)\bigr\rangle_{U,\boldsymbol z}.
\]
This formula identifies the coefficient of each marker by regrouping
the weights on an edge-disjoint decomposition. The marker updates in
Propositions~\ref{prop:torso-marker} and~\ref{prop:partial-branch-records}
supply the required local decision. They distinguish closures with the
same unmarked boundary connectivity and place their ribbon evaluations
in the appropriate sectors. For example, inserting one fixed edge in the
labelled $K_7$ census creates $2{,}220$ Petersen-family certificates.
The four-terminal calculations in Section~\ref{sec:implementation}
compare candidate counts and elapsed times for the two pruning rules. The repeated-cell example of Section~\ref{sec:k6-ear-transfer} gives a closed formula for a full marked sector with arbitrarily many edges.

For generic four-point Temperley--Lieb data, finite return under
each braid slide forces a scalar block. In the Fox model the same
graph rule works for every embedding, uses finite color transport
and yields nonproportional blocks on realizable states. The
remaining generic question is which boundary spaces or admissible
cut moves can support such a coupling while preserving ambient isotopy.

\section*{Data and code availability}
The supplementary source contains \texttt{anc/verify\_examples.py}. It recomputes the square and its sixteen core closures, the labelled $K_6$ and $K_7$ marker counts, and the fixed-edge transition table using the Python standard library. The option \texttt{--cross-check} tests the $K_7$ census by a second method, and \texttt{--coefficients} prints all nonzero classical sector coefficients. 

The companion \texttt{anc/verify\_quantum.py} checks the bracket matrices, the four-port double-slide obstruction, the two cable blocks, and the three-edge deletion--contraction identity using SymPy. The standard-library script \texttt{anc/verify\_spatial\_boundary.py} enumerates all $2^{15}$ states of the decorated $K_6$, checks its sector coefficients and the two link evaluations, compares the four-terminal record update with direct minor-model searches on both hosts, and verifies the two leading $K_7$ coefficients for the Hamiltonian-cycle coloring. It also checks the $K_6$ and $K_7$ minor tests on $K_{3,3,1,1}$
and counts candidate pairs before and after pruning across the four-terminal
labeling classes. The
standalone program \texttt{anc/verify\_records\_partitions.py} checks
the shape multiplicities on four and six terminals and tests the boundary
descriptor count in Section~\ref{sec:four-terminal-models} by direct enumeration.

For the benchmarks, run \path{anc/benchmark_minor.py} with the options
in Section~\ref{sec:record-benchmarks}. The measurements reported there are
in \path{anc/benchmark_observed.csv}. The script \path{anc/verify_transfer_family.py}
checks the two-state transfer in Section~\ref{sec:k6-ear-transfer}. Both
programs use the Python standard library.
The data for Table~\ref{tab:six-block-types} and for the four-terminal
pruning table are in \path{anc/table1_six_boundary_partitions.csv} and
\path{anc/pruning_observed.csv}. To compare Table~\ref{tab:six-block-types} and the two tables in
Section~\ref{sec:implementation} with their CSV rows,
run \path{anc/verify_tables.py} with the PDF path as its argument (Poppler and PyMuPDF
required). The file \path{anc/formulas.txt} gives the linear notation for the
boundary identities and the counting formulas in Table~\ref{tab:six-block-types}
and Section~\ref{sec:complete-marker-census}. To check the formulas and every Greek letter in Definitions~3.1 and~3.4, run
\path{anc/verify_formula_text.py} with the PDF path. The programs \path{anc/verify_finite_orbit.py} and
\path{anc/verify_fox_cover.py} check the projector orbit and the
finite coloring example, respectively.

The moment-curve numbers in Proposition~\ref{stmt:15.3} enter only as inputs to its conditional identity.

\end{document}